\documentclass[10pt]{article}

\RequirePackage{verbatim}
\RequirePackage{hyperref}
\RequirePackage{amsmath,amsfonts,amssymb}

\usepackage{latexsym,amscd,amsthm}
\usepackage[all]{xy}
\usepackage{amstext,amsopn}
\usepackage{graphicx}
\usepackage{amscd}
\usepackage{amssymb}
\usepackage{url}
\usepackage{manfnt}
\usepackage{xcolor}
\usepackage{supertabular}
\usepackage{slashbox}

\usepackage{tikz}

\newtheorem{Thm}{Theorem}[section]
\newtheorem{Prop}[Thm]{Proposition}

\newtheorem{defprop}[Thm]{Definition-Proposition}

\theoremstyle{remark}
\newtheorem{Rem}[Thm]{Remark}

\theoremstyle{definition}

\newtheorem{Def}[Thm]{Definition}
\newtheorem{Exa}[Thm]{Example}

\author{Damien Calaque}
\title{Derived stacks in symplectic geometry}
\date{}

\begin{document}

\maketitle

\tableofcontents

\section*{Introduction}
\addcontentsline{toc}{section}{Introduction}

Our goal in this Chapter is to explain how derived stacks can be useful for ordinary symplectic geometry, with an emphasis on examples coming from classical topological field theories. 
More precisely, we will use classical Chern--Simons theory and moduli spaces of flat $G$-bundles and $G$-local systems as leading examples in our journey. 

\medskip

We will start in the introduction by reviewing various point-of-views on classical Chern--Simons theory and moduli of flat connections. 
In the main body of the Chapter we will try to convince the reader how derived symplectic geometry (after Pantev--To\"en--Vaqui\'e--Vezzosi \cite{PTVV}) somehow reconciles all these different point-of-views. 

\subsection*{Physics: classical Chern--Simons theory}

Let $M$ be a closed oriented $3$-manifold and let $X=Conn_G(M)$ be the space of all $G$-connections on $M$, where $G\subset GL(n,\mathbb{R})$ is a compact simple Lie group. 
\begin{Rem}
We consider connections on arbitrary differentiable principal $G$-bundles here. Nevertheless, as the moduli of differentiable $G$-bundles is discrete, we will deal with connections on the trivial principal $G$-bundle for simplicity in this introduction. 
\end{Rem}
The Chern--Simons action functional is then 
$$
CS(A):=\frac{k}{4\pi}\int_M tr(dA\wedge A+\frac23A\wedge A\wedge A)\,.
$$
Since we will be interested in \textit{classical} (as opposed to \textit{quantum}) Chern--Simons theory then we can safely assume the factor $\frac{k}{4\pi}$ (which is relevant for quantization) is $1$. 

\medskip

The space of classical trajectories of our physical system is given by the space $Crit(CS)$ of critical points of the action functional. 
The Chern--Simons functional has a huge group of symmetries $\mathcal G_M=C^\infty(M,G)$ (gauge symmetries), and we will be essentially interested 
in the ``reduced'' space of trajectories $X_{red}=Crit(CS)/\mathcal G_M$. 

\medskip

If $M$ is now bounded by an oriented closed surface $\Sigma=\partial M$ then we have a restriction map $r:X=Conn_G(M)\to Conn_G(\Sigma)=P$. 
The space $P$ is the phase space of our system and $r$ encodes boundary/initial conditions (or constraints). 
The corresponding reduced phase space $P_{red}$ shall in principle be symplectic and the subspace of admissible boundary/initial conditions 
(i.e.~the image of the induced map $X_{red}\to P_{red}$) shall be Lagrangian (note that we would like to say that $X_{red}$ itself is Lagrangian, but $X_{red}\to P_{red}$ might not be injective, 
i.e.~a gauge equivalence class of classical trajectories might not be uniquely determined by its initial conditions). 

\begin{Rem}\label{rem:infty}
There are infinite dimensional spaces that are involved here. 
We will adopt a functorial approach to differential geometry that will allow us to deal with these infinite dimensional spaces, e.g.~seeing a space as ``its functor of points'': $Conn_G(M)$ is for instance the functor that sends a smooth manifold $X$ to smooth families of $G$-connections on $M$ paramatrized by $X$. We refer to \cite{PIZ} for an approach using \textit{diffeologies}, 
and to \cite{MR} for an approach which uses Dubuc's \textit{$C^\infty$-rings}. The last one has two advantages: it is very close to the modern presentation of algebraic geometry, and it provides a model of synthetic differential geometry (for which we refer to \cite{Kock}).  
\end{Rem}
\begin{Rem}
The reduced spaces, even if they happen to be finite dimensional, may be very singular. 
We will consider them as derived stacks in order to resolve this issue (see \cite{Anel}) . 
\end{Rem}
\begin{Rem}
The problem that $X_{red}\to P_{red}$ is not necessarily injective is fine in the derived setting. 
\end{Rem}

\subsection*{Moduli spaces {\it via} infinite dimensional reduction}

Going back to the case when $M$ is without boundary, one observes that $A\in Conn(M,G)$ is a critical point of the action functional if and only if its curvature $F(A):=dA+A\wedge A$ vanishes. 
In other words, critical points of $CS$ are zeroes of the curvature map $F:X=Conn_G(M)\to \Omega^2(M,\mathfrak g)\subset\Omega^1(M,\mathfrak g)^*$, where the last inclusion is given by 
$\alpha\mapsto \int_Mtr(\alpha\wedge-)$. 
Hence the ``reduced'' space of trajectories $X_{red}=Crit(CS)/\mathcal G_M$ is the quotient $F^{-1}(0)/\mathcal G_M$. 

\medskip

We are tempted to view the curvature map $F$ as a kind of moment map. Up to infinite dimensional issues, this is actually the case if we go down one dimension and consider the curvature map on the phase space $P=Conn_G(\Sigma)$: 
\begin{eqnarray*}
F:P=Conn_G(\Sigma) & \longrightarrow & \Omega^2(\Sigma,\mathfrak g)\subset C^\infty(\Sigma,\mathfrak g)^*\\
A & \longmapsto & dA+ A\wedge A\,.
\end{eqnarray*}
The inclusion $\Omega^2(\Sigma,\mathfrak g)\subset C^\infty(\Sigma,\mathfrak g)^*$ is again given by $\alpha\mapsto\int_{\Sigma}tr(\alpha\wedge-)$. 
Note that $P=Conn_G(\Sigma)$ is a (pre-)symplectic affine space with $2$-form $\omega_P$ given on its tangent $T_A Conn_G(\Sigma)\cong\Omega^1(\Sigma,\mathfrak g)$ by 
$$
\omega_{P,A}(B_1,B_2):=\int_{\Sigma}tr(B_1\wedge B_2)\,.
$$ 
The action of the gauge group $\mathcal G_\Sigma=C^\infty(\Sigma,G)$ on $P$ is Hamiltonian with moment map $F$, 
and the reduced phase space $P_{red}=F^{-1}(0)/\mathcal G_\Sigma$ is obtained {\it via} (infinite dimensional) Hamiltonian reduction. 
\begin{Rem}
Following Remark \ref{rem:infty} the space $P_{red}$ exists as an object in diffeologies. 
As such we can say what a closed $2$-form is, but in order to express what non-degeneracy at a point $[A]\in P_{red}$ means we need to have a nice finite dimensional tangent space $T_{[A]}P_{red}$. 
One can show that there is an open subset of $P_{red}$ where this is the case, and then at such a point $[A]$ the tangent is expressed as follows\footnote{Recall that $P_{red}$ is the quotient of the zero locus of $F$ by $\mathcal G_\Sigma$. }: 
$$
T_{[A]}P_{red}\cong \frac{\ker(dF_A)}{C^\infty(\Sigma,\mathfrak g)}\cong \frac{\ker\big(d_{dR}:\Omega^1(\Sigma,\mathfrak g)\to\Omega^2(\Sigma,\mathfrak g)\big)}{d_{dR}(C^\infty(\Sigma,\mathfrak g)\big)}=H^1(\Sigma,\mathfrak g)\,.
$$
Finally, the reduced pairing $\omega_{P_{red},[A]}$ is non-degenerate on $H^1(\Sigma,\mathfrak g)$, by Poincar\'e duality. 
\end{Rem}
The above is a crucial ingredient in Atiyah--Bott construction of a symplectic structure \cite{AB} on the moduli space. 

\subsection*{Local systems and the quasi-Hamiltonian formalism}

It is a standard fact that flat $G$-bundles up to isomorphisms are exactly local systems up to isomorphisms, i.e.~conjugacy classes of representations of the fundamental group. 
In the case of a closed oriented surface $\Sigma$ the fundamental group $\pi_1(\Sigma)$ admits a presentation with $2g$ generators ($g$ being the genus of $\Sigma$) $a_1,\dots,a_g,b_1,\dots,b_g$ with the sole relation 
\begin{equation}\label{eq:multimoment}
\vec{\prod_{i=1,\dots,g}}(a_i,b_i)=a_0b_0a_0^{-1}b_0^{-1} \cdots a_gb_ga_g^{-1}b_g^{-1}=1\,.
\end{equation}
Therefore the space of $G$-local systems $Loc_G(\Sigma)=Hom_{Grp}\big(\pi_1(\Sigma),G\big)/G$ is of the form $\Phi^{-1}(1)/G$, where $\Phi:G^{2g}\to G$ sends $((a_i,b_i))_{i=1,\dots,n}$ to the l.h.s.~of \eqref{eq:multimoment}. 
One would very much like to see $\Phi$ as a kind of moment map. 
Indeed, Alekseev--Malkin--Meinreinken \cite{AMM} have shown that $\Phi$ is a so-called Lie group valued moment map: in their language $G^{2g}$ is a quasi-Hamiltonian $G$-space and its quasi-Hamiltonian reduction 
(which is genuinely symplectic) is $Loc_G(\Sigma)$. They even interpret $G^{2g}/G$ as $G$-local systems on $\Sigma-D$, where $D$ is a two-dimensional disk, and the map $\Phi^{-1}(1)/G\to G^{2g}/G$ as the restriction map 
$Loc_G(\Sigma)\to Loc_G(\Sigma-D)$, viewing a $G$-local system on $\Sigma$ as a $G$-local system on $\Sigma-D$ with trivial holonomy around $\partial D$. 

\medskip

Observe that $G$ itself is the quotient $C^\infty(S^1,G)/C^\infty_*(S^1,G)$, where $C^\infty_*(S^1,G)=L_*(G)$ denotes based loops (sending $1$ to $1$). 
This allows to identify $G^{2g}$ as the moduli of flat connections on $\Sigma-D$ up to gauge equivalences fixing a point in the boundary. 
Finally, the space $FlatConn_G(\Sigma-D)$ of all flat connections on $\Sigma-D$ is equipped with a restriction map 
$$
\mu:FlatConn_G(\Sigma-D)\to FlatConn_G(S^1)=\Omega^1(S^1,\mathfrak g)\subset C^\infty(S^1,\mathfrak g)^*\,,
$$
which is a moment map for the action of $L(G)=C^\infty(S^1,G)$. 
Using the fact that $L(G)/L_*(G)=G$, together with the holonomy map $\Omega^1(S^1,\mathfrak g)\to G$, Alekseev--Malkin--Meinreinken \cite{AMM} show that there is a correspondence between 
Hamiltonian $L(G)$-spaces and quasi-Hamiltonian $G$-spaces, and that there are symplectomorphisms: 
$$
P_{red}\,\tilde\longleftarrow\,\mu^{-1}(0)/L(G)\,\tilde\longrightarrow\,\Phi^{-1}(1)/G\,,
$$
where the left arrow sends a flat connection on $\Sigma-D$ that is zero on the boundary to its extension by zero, and the right arrow sends a flat connection to its holonomy representation. 

\medskip

The work of Alekseev--Malkin--Meinreinken provides a finite dimensional construction of the symplectic structure on the moduli space of flat connections, for which the topological invariance of the reduced phase space is furthermore transparent. 

\subsection*{Deformation theoretic approach}

There is a clean way of showing that there is a non-degenerate pairing on the tangent to the moduli space, viewed as a functor of points, of either flat $G$-bundles or $G$-local systems. We have that 
$$
T_{[\mathcal P,\nabla]}Flat_G(\Sigma)=H^1_{dR}\big(\Sigma,(\mathcal P\times_G\mathfrak g,\nabla)\big)
$$
is the first de Rham cohomology of $\Sigma$ with values in the associated flat vector bundle $\mathcal P\times_G\mathfrak g$. 
and
$$
T_{[\mathcal L]}Loc_G(\Sigma)=H^1_{Betti}(\Sigma,\mathcal L\times_G\mathfrak g)\,.
$$
is the first singular/Betti cohomology of $\Sigma$ with values in the associated linear local system $\mathcal L\times_G\mathfrak g$. 

In both cases, using the trace pairing together with the cup-product we get a non-degenerate skew-symmetric form with values in $H^2_{Betti}(\Sigma,\mathbb{R})\cong H^2_{dR}(\Sigma,\mathbb{R})\cong\mathbb{R}$. 
On can easily show that it varies smoothly with respect to $[\mathcal P,\nabla]$ or $[\mathcal L]$. But proving that it is closed is not easy (and basically reduces to make use of one of the previous approaches). 

\subsection*{Back to physics: local critical loci}

Going along the lines of \cite{JS,KL} one can show that the moduli space of flat $G$-connections on a compact oriented $3$ manifold is locally the critical locus of a function on a finite dimensional differentiable space\footnote{In \cite{JS,KL} this is done for holomorphic Chern--Simons theory. I don't know a reference where this is done for Chern--Simons theory, but it shall be easier than its holomorphic analog. }. In order to see this, one has to observe that near a given flat $G$-connection $\nabla$ one can put several constraints (like Laplacian eigenvalue constraints) that allows to identify a neighborhood of $\nabla$ as a critical locus of the Chern--Simons functional restricted to a finite dimensional constrained subspace. 

This allows in particular to get rid of infinite dimensional complications, at the price of loosing the existence of a global action functional. Nevertheless, one can still use this local structure in order to compute interesting invariants. 

\subsection*{Unifying all these approaches}

Our goal in this Chapter is to present a context in which the following questions can be addressed in a meaningful way and get natural answers: 
\begin{itemize}
\item[\textbf{[a]}] how to relate these different approaches? 
More precisely, how to put a natural symplectic structure on our favorite moduli so that all these different approaches somehow appear as {\it different ways of computing} this symplectic structure?
\item[\textbf{[b]}] whenever $\partial M=\emptyset$, what does it mean for $X_{red}\cong Loc_G(M)$ to be Lagrangian in $P_{red}=*$?\footnote{The analogous situation in classical mechanics would be to look at the space of periodic trajectories of a classical mechanical system, which is also a typical example where classical trajectories are not determined by initial/boundary conditions. }
\item[\textbf{[c]}] how to get rid of smoothness issues? Namely, in the Introduction we have ``carefully'' avoided to be precise about the fact that most statements are true only when one is over some smooth locus in moduli spaces.  
\item[\textbf{[d]}] how to get rid of the problem that $r:X_{red}\to P_{red}$ may not be an inclusion (for instance, when $P_{red}$ is a point). 
\item[\textbf{[e]}] to what extent (and how) one can recover an action functional from $X_{red}$? 
\item[\textbf{[f]}] how to express the fact that Chern--Simons theory is a topological field theory?
\end{itemize}

\subsection*{Our geometric context}

We refer to Mathieu Anel's contribution \cite{Anel} for an introduction to the ideas of derived geometry. 
We will work in the $C^\infty$ setting, using for instance the theory of derived differentiable stacks from \cite{Wa} (see also \cite{DAG-V}), which uses $C^\infty$-rings. 

Smoothness issues appearing in question \textbf{[c]} will be resolved by this use of derived differentiable stacks. 
There is a large class of derived differentiable stacks, the Artin ones, on which all the calculations we want to perform make sense and that is stable both by pull-backs and by groupoid action quotients. 
Surprisingly enough, the other questions will get very natural answers once one has set up a suitable theory of symplectic structures on Artin stacks. 

\subsection*{Acknowledgements}

First of all I warmly thank Mathieu Anel for the many discussions we had about derived stacks in symplectic geometry (and plenty other topics). I also thank him, as well as Gabriel Catren, for their very constructive and helpful comments on a preliminary version of this Chapter. I then owe a debt of gratitute to Tony Pantev, Bertrand To\"en, Michel Vaqui\'e and Gabriele Vezzosi. Finally, numerous discussions with Pavel Safronov have been very useful. 

The idea of having classical Chern--Simons theory as a common thread came up when I was preparing a colloquium talk\footnote{\textit{Derived symplectic geometry and classical Chern--Simons theory}, \url{http://www.perimeterinstitute.ca/fr/videos/derived-symplectic-geometry-and-classical-chern-simons-theory}. } for a conference on \textit{Deformation quantization of shifted Poisson structures} at Perimeter Institute. 

This work has been partly supported by the \textit{Institut Universitaire de France} and by the \textit{Agence Nationale de la Recherche} (through the project ``Structures sup\'erieures en Alg\`ebre et Topologie''). 


\section{Symplectic structures in the derived setting}\label{section1}

Recall that a \emph{symplectic structure} on a smooth manifold $M$ is a $2$-form $\omega\in\Omega^2(M)$ that is 
\begin{itemize}
\item closed for the de Rham differential: $d_{dR}\omega=0$. 
\item non-degenerate: the contraction map $\omega^b:TX\to T^*X,v\mapsto\iota_v\omega$ is a bundle isomorphism. 
\end{itemize}
A submanifold $L\subset M$ is \emph{Lagrangian} if it is
\begin{itemize}
\item isotropic: $\omega_{|L}$ vanishes on $TL$. 
\item non-degenerate: the induced map $TL\to N^*L$, where $N^*L$ is the conormal bundle of $L$, is a bundle isomorphism.  
\end{itemize}

\subsection{Closed forms: structure {\it versus} property}

Working in the derived setting, one shall have a notion of differential form that is homotopy invariant in the sense that it does not depend on the explicit presentation of our derived stack. More precisely, Artin stacks have a cotangent \textit{complex}, so that forms will have a cohomological degree. Different presentations of the same stack will lead to quasi-isomorphic cotangent complexes. Therefore if we get a form that is de Rham closed in the naive sense for a given presentation, it might only be de Rham closed \textit{up to a cocycle} for another presentation\footnote{This is rather usual that algebraic identites are not strictly preserved under quasi-isomorphisms.}. 
This leads to the idea of considering forms that are closed \textit{up to homotopy} rather than in the strict (naive) sense. 
We will try to implement this idea and show that it naturally follows from the synthetic and functorial approach to differential forms. 

For the first reading, the reader who is only interested in a concrete and intuitive definition of these closed forms up to homotopy can directly go to the paragraph \ref{sssec1.5} dealing with examples, and consider that a \textit{closed $p$-form of degree $n$} on a derived stack $X$ is a series $(\alpha_p,\alpha_{p+1},\dots,)$ where 
\begin{itemize}
\item $\alpha_i$ is an $i$-form of total degree $p+n$, where the total degree is the sum of the cohomological degree and the form degree\footnote{The form degree will be named \textit{weight}, and the total degree will simply be named \textit{degree}. }. 
\item $d_{dR}(\alpha_i)=\partial(\alpha_{i+1})$, where $d_{dR}$ is the de Rham differential and $\partial$ is the \textit{internal} differential (that comes from that we are doing derived geometry), and with the convention that  $\partial(\alpha_p)=0$. 
\end{itemize}
The leading term $\alpha_p$ is called the underlying $p$-form of degree $n$. 

\subsubsection{The de Rham complex from the synthetic point-of-view}\label{sssec1.1}

For a genuine smooth manifold $X$, one can consider the differentiable space $X^{\Delta^k_{inf}}$ of maps from the infinitesimal $k$-simplex $\Delta^k_{inf}$ to $X$. Recall that the infinitesimal $k$-simplex $\Delta^k_{inf}$ is the space of $k+1$-tuples of points $(x_0,\dots, x_k)$ in $\mathbb{R}^k$ that are pairwise close at order $1$; in other words, its $C^\infty$-ring of functions is the quotient of $C^\infty(\mathbb{R}^{k+1})$ by the relations
$$
\sum_{i=0}^kx_i=1\quad\textrm{and}\quad (x_i-x_j)^2=0~\textrm{if}~i\neq j\,.
$$
For instance, the infinitesimal $1$-simplex $\Delta^1_{inf}$ is a so-called ``fat point'', having function $C^\infty$-ring the quotient of $C^\infty(\mathbb{R})$ by $t^2=0$ (here $t$ is $x_0-x_1$), and thus $X^{\Delta^1_{inf}}$ is the tangent space $TX$ to $X$. 

The collection $(\Delta^k_{inf})_{k\geq0}$ is cosimplicial, hence the collection $(X^{\Delta^k_{inf}})_{k\geq0}$ forms a simplicial object in differentiable spaces and thus, $C^\infty(X^{\Delta^k_{inf}})_{k\geq0}$ is a cosimplicial vector space. 
It is a fact (see \cite[\S I.18]{Kockbook}) that the usual differential $k$-forms $\Omega^k(X)$ may be characterized as the joint kernel of the degeneracy maps in $C^\infty(X^{\Delta^k_{inf}})$. 
Hence the de Rham complex $DR(X)$ of $X$ is identified with the normalized complex\footnote{Recall that the normalized complex of a cosimplicial vector space $(V^k)_{k\geq0}$ is constructed as follows: in degree $k$ one has the joint kernel of degeneracies $V^k\to V^{k-1}$, and the differential is the alternating sum of faces $V^k\to V^{k+1}$. } of the cosimplicial vector space $C^\infty(X^{\Delta^k_{inf}})_{k\geq0}$. 

\medskip

\noindent Closed $p$-forms on $X$ can then be characterized in the following equivalent ways: 
\begin{itemize}
\item[\textbf{(c1)}] as $p$-cocycles in $DR(X)$. 
\item[\textbf{(c2)}] as morphisms of complexes from $\mathbb{R}[-p]$ to $DR(X)$.\footnote{For a cochain complex $C$ and an integer $k$, its shift by $k$, denoted $C[k]$, is another cochain complex whose degree $m$ cochains are the degree $m+k$ cochains of $C$. For instance, $\mathbb{R}[-p]$ is the cochain complex that only consists of a one dimensional space of degree $p$ cochains. } 
\item[\textbf{(c3)}] as morphisms of cosimplicial vector spaces from $\mathbb{R}(-p):=K(\mathbb{R}[-p])$ to $C^\infty(X^{\Delta^k_{inf}})_{k\geq0}$, where $K$ is inverse to the normalized complex functor. 
\end{itemize}
Note that the cosimplicial vector space $\mathbb{R}(-p)$ admits the following very simple description: 
in cosimplicial degree $k<p$ it is zero, in cosimplicial degree $p$ it is $\mathbb{R}$, and it is freely generated by faces in higher cosimplicial degree. 

\medskip

\noindent Observe that (non-necessarily closed) $p$-forms can be characterized in similar terms: 
\begin{itemize}
\item[\textbf{(nc1)}] as $p$-cochains in $DR(X)$. 
\item[\textbf{(nc2)}] as morphisms of graded vector spaces from $\mathbb{R}[-p]$ to $DR(X)^\sharp$, where $(-)^\sharp$ stands for the underlying graded vector space functor. 
\item[\textbf{(nc3)}] as morphisms of $p$-truncated cosimplicial vector spaces from $\mathbb{R}(-p)_{\leq p}$ to $C^\infty(X^{\Delta^k_{inf}})_{0\leq k\leq p}$. 
\end{itemize}

\noindent We want to generalize the above three definitions for derived differentiable stacks, from the third to the first. 
\begin{Rem}
The first definition is obviously the simplest and most concrete of the three. 
Its generalization appears to be rather intuitive and explicit, but proving general results about closed forms on derived stacks with this definition can be rather complicated. The other definitions as well as their generalizations are rather abstract, but they appear to be more convenient for proving general statements about closed forms on derived stacks. 

As we will see, the third definition has a straightforward generalization. The generalization of the second one is less straightforward, but we give it for several reasons: it is rather compact, it is the one that nowadays is used in most reference (see e.g.~\cite{CPTVV}), and it explains the relation between the third and the first definitions. 
\end{Rem}

\subsubsection{Closed forms on derived stacks: third definition}

If $X$ is a derived differentiable stack, then we have a derived stack 
\[
X^{\Delta^k_{inf}}=\mathbf{Map}(\Delta^k_{inf},X)
\]
of maps from $\Delta^k_{inf}$ to $X$, whose $Y$-points are maps $\Delta^k_{inf}\times Y\to X$. Recall that derived global functions on a derived stack form a differential graded algebra, and thus, in particular, a cochain complex. 
Hence $C^\infty(X^{\Delta^k_{inf}})$ is a cochain complex, and $C^\infty(X^{\Delta^k_{inf}})_{k\geq0}$ is a cosimplicial cochain complex (i.e.~a cosimplicial object in cochain complexes). 
\[
\xymatrix{
\cdots\ar[r]^-{\partial} & C^\infty(X^{\Delta^2_{inf}})^i \ar@<-2ex>[d] \ar@<2ex>[d] \ar[r]^-{\partial} & 
C^\infty(X^{\Delta^2_{inf}})^{i+1} \ar@<-2ex>[d] \ar@<2ex>[d] \ar[r]^-{\partial} & \cdots \\
\cdots \ar[r]^-{\partial} & 
C^\infty(X^{\Delta^1_{inf}})^i \ar@<-4ex>[u] \ar[u] \ar@<4ex>[u] \ar[d]\ar[r]^-{\partial} & 
C^\infty(X^{\Delta^1_{inf}})^{i+1} \ar@<-4ex>[u] \ar[u] \ar@<4ex>[u]  \ar[d] \ar[r]^-{\partial} & \cdots \\
\cdots \ar[r]^-{\partial} & 
C^\infty(X)^i \ar@<-2ex>[u] \ar@<2ex>[u] \ar[r]^-{\partial} & 
C^\infty(X)^{i+1}\ar@<-2ex>[u] \ar@<2ex>[u] \ar[r]^-{\partial} & \cdots
}
\]
This suggests the following generalization of characterization \textbf{(c3)} for closed $p$-forms: 
\begin{Def}
A closed $p$-form of degree $n$ is a morphism of cosimplicial cochain complexes from $\mathbb{R}(-p)$ to $\big(C^\infty(X^{\Delta^k_{inf}})[n]\big)_{k\geq0}$. 
One similarly defines a (non-necessarily closed) $p$-form of degree $n$ as a morphism of $p$-truncated cosimplicial cochain complexes from $\mathbb{R}(-p)_{\leq p}$ to $\big(C^\infty(X^{\Delta^k_{inf}})[n]\big)_{0\leq k\leq p}$. 
\end{Def}
\noindent Note that we have a space\footnote{Here we mean a space in the sense of homotopy theory, i.e.~a homotopy type. } $\mathcal A^{p,cl}(X,n)$ of such closed $p$-forms of degree $n$. For instance, a path between two of these is a homotopy between the corresponding maps of cosimplicial cochain complexes\footnote{More 
precisely, we have a class of weak equivalences in the category of cosimplicial cochain complexes: the levelwise quasi-isomorphisms. Localizing at weak equivalences we therefore have an $\infty$-category of cosimplicial cochain complexes, 
and thus we have spaces of morphisms. These spaces of morphisms can be computed using an explicit model structure. }. 
Similarly, we write $\mathcal A^{p}(X,n)$ for the space of $p$-forms of degree $n$. 
There is an obvious map $\mathcal A^{p,cl}(X,n)\to \mathcal A^{p}(X,n)$ which is not necessarily a subspace map. 
This means that for derived stacks (as opposed to genuine manifolds) being closed is a structure rather than just a property. 
We will make this more transparent below. 

\subsubsection{Closed forms and graded mixed complexes: second definition}

We have seen that, when $X$ is a derived differentiable stack, each $C^\infty(X^{\Delta^k_{inf}})$ is already a cochain complex. 
We can still consider the normalized complex $DR(X)$, which has the richer structure of a \textit{graded mixed complex}\footnote{The notion of a \textit{mixed complex} (a term coined by Kassel \cite{Kassel}) goes back to Dwyer and Kan's study of Connes' cyclic modules\cite{DK}, who were calling it a \textit{duchain complex}. \textit{Graded mixed complexes} have been introduced in \cite{PTVV}, but we would like to warn the reader that in this Chapter we adopt the grading convention from \cite{CPTVV}. }: 
\begin{itemize}
\item it is a complex, as such it carries a grading that we will refer to as the \textit{cohomological grading}. 
\item it carries an auxiliary grading, the \textit{weight}, that is reminiscent from the cosimplicial degree: element in 
$C^\infty(X^{\Delta^k_{inf}})$ (i.e.~$k$-forms) have weight $k$. 
\item its differential splits into a sum $\partial+\epsilon$: 
\begin{itemize}
\item[a)] $\partial$ is the part which comes from the differential of $C^\infty(X^{\Delta^k_{inf}})$, that we call the \textit{internal differential}. 
It is of zero weight. 
\item[b)] $\epsilon$ is the part which comes from the alternating sum of the face maps, that we call in general the \textit{mixed differential}, and the \textit{de Rham differential} $d_{dR}$ in this specific example. 
It has weight one. 
\end{itemize}
\end{itemize}
If one denotes by $DR(X)_{(p)}$ the weight $p$ part of $DR(X)$, then the graded mixed complex $DR(X)$ can be visualized as follows:\footnote{Observe that $DR(X)_{(0)}=C^\infty(X)$. } 
\[
\xymatrix{
\cdots \ar[r]^-{\partial} & 
DR(X)_{(2)}^i \ar[r]^-{\partial} & 
DR(X)_{(2)}^{i+1} \ar[r]^-{\partial} & \cdots \\
\cdots \ar[r]^-{\partial} \ar[ur]^-{d_{dR}} & 
DR(X)_{(1)}^i \ar[ur]^-{d_{dR}} \ar[r]^-{\partial} & 
DR(X)_{(1)}^{i+1}\ar[ur]^-{d_{dR}} \ar[r]^-{\partial} & \cdots \\
\cdots \ar[r]^-{\partial} \ar[ur]^-{d_{dR}} & 
DR(X)_{(0)}^i \ar[ur]^-{d_{dR}} \ar[r]^-{\partial} & 
DR(X)_{(0)}^{i+1}\ar[ur]^-{d_{dR}} \ar[r]^-{\partial} & \cdots
}
\]
\noindent The normalized complex functor $N$ induces an equivalence of categories between cosimplicial complexes and graded mixed complexes sitting in non-negative weight, 
which induces an equivalence between their $\infty$-categorical localizations\footnote{Weak equivalences are levelwise quasi-isomorphisms in both cases. }
\begin{Rem}
Note that $N\big(\mathbb{R}(-p)\big)=\mathbb{R}[-p](-p)$ is just $\mathbb{R}$ viewed as a graded mixed complex concentrated in weight $p$ and cohomological degree $p$. 
\end{Rem}
One therefore has the following generalization of characterization \textbf{(c2)} for closed $p$-forms, which is precisely the one appearing in \cite{CPTVV} (and most recent references):
\begin{defprop}
The space $\mathcal A^{p,cl}(X,n)$ of closed $p$-forms of degree $n$ is the space of maps from $N\big(\mathbb{R}(-p)\big)$ to $DR(X)[n]$. 
The space $\mathcal A^{p}(X,n)$ of $p$-forms of degree $n$ is the space of maps from $N\big(\mathbb{R}(-p)\big)^\sharp$ to $DR(X)[n]^\sharp$, where $(-)^\sharp$ stands for the underlying graded complex functor (forgetting the mixed differential). The map $\mathcal A^{p,cl}(X,n)\to \mathcal A^{p}(X,n)$ comes from applying $(-)^\sharp$. 
\end{defprop}
\begin{Rem}
In the non-derived setting (see paragraph \ref{sssec1.1}) the graded mixed complex $DR(X)$ is diagonal (in the sense that there is no internal differential and its weight $p$ part is concentrated in cohomological degree $p$). 
\end{Rem}
\begin{Rem}\label{bicpx}
Graded mixed complexes are very much related to filtered complexes and spectral sequences. 
For instance, with every bi-complex $(C^{\bullet,\bullet},d_h,d_v)$ one can associate two graded mixed complexes: 
\begin{itemize}
\item[(1)] the cohomological grading is the total degree and the weight is the horizontal degree; the internal differential is $d_v$ and the mixed differential is $d_h$. 
\item[(2)] the cohomological grading is the total degree and the weight is the vertical degree; the internal differential is $d_h$ and the mixed differential is $d_v$. 
\end{itemize}
Conversely, one can build two bi-complexes from a graded mixed complex $(D,\partial,\epsilon)$ by saying that: 
\begin{itemize}
\item[(1)] an element of weight $p$ and degree $d$ has bidegree $(p,d-p)$, $d_h=\epsilon$ and $d_v=\partial$. 
\item[(2)] an element of weight $p$ and degree $d$ has bidegree $(d-p,p)$, $d_h=\partial$, and $d_v=\epsilon$. 
\end{itemize}
For instance, in our example of the de Rham graded mixed complex $DR(X)$, letting $\Omega^p(X):=DR(X)_{(p)}[p]$ we get a bicomplex 
\[
\xymatrix{
\cdots \ar[r]^-{\partial} & 
\Omega^2(X)^i \ar[r]^-{\partial} & 
\Omega^2(X)^{i+1} \ar[r]^-{\partial} & \cdots \\
\cdots \ar[r]^-{\partial} & 
\Omega^1(X)^i \ar[u]^-{d_{dR}} \ar[r]^-{\partial} & 
\Omega^1(X)^{i+1}\ar[u]^-{d_{dR}} \ar[r]^-{\partial} & \cdots \\
\cdots \ar[r]^-{\partial} & 
C^\infty(X)^i \ar[u]^-{d_{dR}} \ar[r]^-{\partial} & 
C^\infty(X)^{i+1}\ar[u]^-{d_{dR}} \ar[r]^-{\partial} & \cdots
}
\]
\end{Rem}

\subsubsection{Closed forms as cocycles: first definition}\label{sssec1.4}

We would like to have a more concrete definition of (closed) forms, and have an explicit model for the space of morphisms of graded mixed complexes from $N\big(\mathbb{R}(-p)\big)$ to $DR(X)$ (or any other graded mixed complex). This is a rather standard problem people encounter in homological and homotopical algebra: one needs to replace $N\big(\mathbb{R}(-p)\big)$ with an equivalent graded mixed complex (a \textit{resolution}) $R_p$ such that mapping out of it is better behaved\footnote{For the expert reader, $R_p$ will be an explicit cofibrant replacement of $N\big(\mathbb{R}(-p)\big)$ in the projective model structure on graded mixed complexes.}. Following \cite{CPTVV} we have the following explicit nice replacement $R_p$ for $N\big(\mathbb{R}(-p)\big)$: it is the linear span of $\{x_i,y_i\}_{i\geq p}$ with: 
\begin{itemize}
\item $x_i$'s having degree $p$ and weight $i$. 
\item $y_i$'s having degree $p+1$ and weight $i+1$. 
\item $\partial(x_i)=y_{i-1}$ and $\epsilon(x_i)=y_i$ (convention: $y_{p-1}=0$). 
\end{itemize}
In particular one sees that $\partial(x_p)=0$ and $\epsilon(x_i)=\partial(x_{i+1})$ for every $i\geq p$. 
\begin{center}
\begin{tabular}{c|c|c|c|ccc|c|cc}
\textrm{degree $\backslash$ weight} & $p$   & $p+1$     & $p+2$     && $\cdots$ && $i$       & $i+1$     \\
\hline
$p+1$                         &     & $y_{p}$   & $y_{p+1}$ && $\cdots$ && $y_{i-1}$ & $y_i$     \\
\hline
$p$                           & $x_p$ & $x_{p+1}$ & $x_{p+2}$ && $\cdots$ && $x_i$     & $x_{i+1}$ \\ 
\end{tabular}
\end{center}
\begin{Rem}
If one accepts that free objects are nice, then one is lead to construct $R_p$ in the following way: 
\begin{itemize}
\item introduce $y_p$ as freely generated by $\epsilon$ from $x_p$: $\epsilon(x_p)=y_p$. 
\item the new graded mixed complex is no longer equivalent to the original one, hence we introduce $x_{p+1}$ in order to ``kill'' the new cohomology class corresponding to $y_p$: $\partial(x_{p+1})=y_p$. 
\item iterate the process. 
\end{itemize}
\end{Rem}
Thus an element in the space $\mathcal A^{p,cl}(X,n)$ can be represented by a genuine morphism of graded mixed complexes $R_p\to DR(X)$. 
It consists in a collection $(\alpha_p,\alpha_{p+1},\dots)$, where $\alpha_i$ has weight $i$ and degree $p+n$, 
such that $\partial(\alpha_p)=0$ and $d_{dR}(\alpha_i)=\partial(\alpha_{i+1})$ for every $i\geq p$. 
Its underlying $p$-form is the leading term $\alpha_p$. 

\medskip

\noindent Let $DR(X)_{(\geq p)}:=\prod_{i\geq p}DR(X)_{(i)}$ be the completed weight $\geq p$ part of $DR(X)$. 
We have the following generalization of characterization \textbf{(c1)} for closed $p$-forms:
\begin{defprop}
The space $\mathcal A^{p,cl}(X,n)$ of closed $p$-forms of degree $n$ is the space of $(p+n)$-cocycles in $DR(X)_{(\geq p)}$ for the total differential $\partial+d_{dR}$. 
The space $\mathcal A^{p}(X,n)$ of $p$-forms of degree $n$ is the space of $(p+n)$-cocycles in $DR(X)_{(p)}$ for the differential $\partial$. 
The map $\mathcal A^{p,cl}(X,n)\to \mathcal A^{p}(X,n)$ is given by extracting the leading term. 
\end{defprop}

\noindent This definition is roughly the one that originally appears in \cite[\S1.2]{PTVV}. 

\subsubsection{Examples}\label{sssec1.5}

Let us summarize what is known about the de Rham complex: 
\begin{itemize}
\item if $X$ is a genuine smooth manifold, then $DR(X)$ is the diagonal graded mixed complex associated with the usual de Rham complex 
$(\Omega^*(X),d_{dR})$. 
\item for a {\it geometric} derived stack $X$, the weight $p$ part $DR(X)_{(p)}$ can be computed in terms of the cotangent complex $\mathbb{L}_X$: $DR(X)_{(p)}\cong\Gamma\big(S^p(\mathbb{L}_X[-1])\big)$ is the $p$-th symmetric power of a shift of the cotangent complex\footnote{The cotangent complex always exists for a derived geometric stack. When it exists, the cotangent complex $\mathbb{L}_X$ of a derived stack $X$ is a far reaching generalization of the cotangent bundle of a manifold, that encodes the infinitesimal structure of the derived stack $X$. We view it as a $\mathcal O_X$-module, where $\mathcal O_X$ is the sheaf of functions on $X$ (a sheaf of differential graded commutative algebras). Tensor, symmetric and skew-symmetric powers of $\mathcal O_X$-modules are understood over $\mathcal O_X$. }. 
In other words, $DR(X)^\sharp\cong\Gamma\big(S(\mathbb{L}_X[-1](-1))\big)$ as graded complexes. 
\item the \textit{de Rham functor} $DR$ satisfies smooth descent: if $X\cong colim_{[n]\in\Delta^{op}} X_n$ is a geometric derived stack that is presented by a nice enough simplicial diagram $X_\bullet$ in stacks then $DR(X)\cong lim_{[n]\in\Delta}DR(X_n)$. Here are examples of nice enough simplicial diagrams of stacks: 
\begin{itemize}
\item $X_\bullet$ is the nerve of a Lie groupoid $G_1\rightrightarrows G_0$: in simplicial degree $n$, 
$$
X_n=\underbrace{G_1\times_{G_0}\cdots\times_{G_0}G_1}_{n\textrm{ times}}\,.
$$
Recall that a Lie groupoid is a groupoid object in manifolds such that the source and target maps $G_1\to G_0$ are submersions (submersions are also called smooth morphisms in algebraic geometry). 
\item $X_\bullet$ is a submersive Segal groupoid in the sense of To\"en-Vezzosi \cite[\S1.3]{TV} (nerves of Lie groupoids that we just mentionned are examples of these). 
\item $X_\bullet$ is a Lie $n$-groupoid after \cite[Definitin 1.2]{Zhu}, i.e.~a simplicial manifold for which horn maps $h_{q,k}$ are surjective submersions for $q\geq1$ and diffeomorphisms for $q\geq n$. 
\end{itemize}
\end{itemize}
This in particular tells us that, if $X=[G_0/G_1]$ is the quotient stack of a genuine Lie groupoid $G_1\rightrightarrows G_0$, then 
$$
DR(X)\cong lim_{[n]\in\Delta}DR(\underbrace{G_1\times_{G_0}\cdots\times_{G_0}G_1}_{n\textrm{ times}})
$$
is equivalent to one of the two graded mixed complexes associated with the following first quadrant bi-complex (according to Remark \ref{bicpx}) from \cite[(2)]{Xu}: 
$$
\xymatrix{
C^\infty(G_1\times_{G_0}\times G_1) \ar[r]^{d_{dR}} & \Omega^1(G_1\times_{G_0}\times G_1) \ar[r] & \cdots                    & \cdots \\
C^\infty(G_1)\ar[u]^{\partial}\ar[r]^{d_{dR}}  		 & \Omega^1(G_1) \ar[u]^{\partial}\ar[r]^{d_{dR}}       & \Omega^2(G_1)\ar[u]\ar[r] & \cdots \\
C^\infty(G_0)\ar[u]^{\partial}\ar[r]^{d_{dR}}			 & \Omega^1(G_0) \ar[u]^{\partial}\ar[r]^{d_{dR}}       & \Omega^2(G_0)\ar[u]^{\partial}\ar[r]^{d_{dR}} & \Omega^3(G_0) \ar[u]
}
$$
Here the horizontal/mixed differential is the usual de Rham differential $d_{dR}$ and the vertical/internal differential is the alternating sum $\partial$ of pull-backs of forms along coface maps of the nerve of the groupoid $G_1\rightrightarrows  G_0$. The degree is the total degree while the weight is the degree of forms. 
\begin{Exa}[Pre-quasi-symplectic groupoids]\label{ex1.6}
Let us make explicit what it means to have a closed $2$-form $\omega$ of degree $1$ on a quotient stack $X=[G_0/G_1]$ as above. 
It is a cocycle of degree $3$ in the total complex that is of weight at least $2$ (the weight being here the genuine degree of forms). 
It is thus an element $\omega=\omega_0+\omega_1\in \Omega^2(G_1)\oplus \Omega^3(G_0)$ such that $\partial\omega_0=0$, $d_{dR}\omega_0=\partial\omega_1$ and $d_{dR}\omega_1=0$. 
This is exactly the notion of pre-quasi-symplectic groupoid from \cite[Definition 2.1]{Xu} (also called a twisted presymplectic groupoid in \cite[\S2.1]{BCWZ}). 
\end{Exa}

If $G_1\rightrightarrows G_0$ is an action groupoid $G\times M\rightrightarrows M$, with $G$ a compact Lie group, $DR(X)$ can then be described as follows: 
\begin{itemize}
\item $DR(X)=\big(S^\star(\mathfrak{g}^*)\otimes\Omega^\bullet(M)\big)^G$, with cohomological grading $2\star+\bullet$ and weight $\star+\bullet$, which we may view as polynomial functions on $\mathfrak g$ with values in differential forms on $M$. 
\item the internal differential is defined as $\partial(\omega)(x):=\iota_{\vec{x}}(\omega(x))$, where $\vec{x}$ is the fundamental vector field associated with $x\in\mathfrak{g}$.  
\item the mixed differential is the de Rham differential on $M$: $\epsilon(\omega)(x):=(d_{dR}\omega)(x)$. 
\end{itemize}
Its total complex is the Cartan model for equivariant cohomology and the weight produces the filtration from \cite[Equation (21)]{Mein} which gives rise to the algebraic counterpart of the Leray spectral sequence. 
\begin{Rem}
The above is consistent with the fact that $DR(X)^\sharp\cong\Gamma\big(S(\mathbb{L}_X[-1](-1))\big)$. Namely, the cotangent complex of a quotient stack $X=[M/G]$ is the two-term complex of $G$-equivariant vector bundles
\begin{eqnarray*}
0 & & 1 \\
T^*_M & \longrightarrow & \mathfrak{g}^*\times M
\end{eqnarray*}
on $M$, where the differential sends a co-vector $\alpha\in T^*_mM$ to the linear map $x\mapsto \alpha(\vec{x}_m)$, and $\vec{x}$ is the fundamental vector field associated with $x\in\mathfrak{g}$. 
When $G$ is compact there is no group cohomology, so that derived global sections over $X$ are just $G$-invariant global sections over $M$, and we are done.  
If $G$ is not compact, then the story is more complicated as it involves non-trivial group cohomology, but we still have a map of graded mixed complexes 
$$
\big(S^\star(\mathfrak{g}^*)\otimes\Omega^\bullet(M)\big)^G\to DR(X)\,.
$$
\end{Rem}
\begin{Exa}[Closed $2$-forms of degree $2$ on $BG$]\label{ex-bg}
Closed $2$-forms of degree $2$ on $BG=[*/G]$ are exactly $G$-invariant symmetric bilinear forms, i.e.~elements in $S^2(\mathfrak g^*)^G$. 
\end{Exa}

As we have seen, if $X\cong colim_{[n]\in\Delta^{op}} X_n$ is a geometric stack that is presented by a Lie $n$-groupoid $X_\bullet$, then we again have that $DR(X)\cong lim_{[n]\in\Delta}DR(X_n)$. In particular, we see that (normalized) multiplicative $p$-forms on $X_q$ (see e.g.~\cite[Definition 2.3]{MT}) are exactly $p$-forms of degree $p-q$ on $X$. A closed $p$-form of degree $p-q$ on $X$ would then be a (normalized) element 
$$
\omega_0+\cdots+\omega_q\in\Omega^p(X_q)\oplus\cdots\oplus\Omega^{p+q}(X_0)
$$
that is closed under the total differential $d_{dR}+\partial$ (where $\partial$ is the internal differential, and is again given explicitly by the alternating sum of pull-backs of forms along coface maps). 

\subsection{Non-degeneracy: shifted symplectic structures}

\subsubsection{Symplectic linear algebra in the $\infty$-categorical setting}

In this paragraph we follow closely the presentation  from \cite{Cal17}, where one works within a general ambient stable symmetric monoidal $\infty$-category $(\mathcal C,\otimes,\mathbf{1})$. Here we will restrict our attention to the following two examples of such $\infty$-categories: 
\begin{itemize}
\item cochain complexes, with equivalences being quasi-isomorphisms, monoidal product the usual tenor product of cochain complexes, and monoidal unit being $\mathbb{R}$. 
\item a sheafified verison of the above over a given stack $X$, denoted $QCoh(X)$. It is a bit larger than the category of complexes of vector bundles over $X$ (which does not necessarily admits fiber products). Its monoidal product is the (derived) tensor product of sheaves of $\mathcal O_X$-modules, and its monoidal unit is $\mathcal O_X$. 
\end{itemize}
As a matter of notation, we recall that we write $?[1]$ for the degree shift functor. 
We define an \textit{$n$-shifted pre-symplectic} object as a pair $(V,\omega)$ where $V$ is an object of $\mathcal C$ and $\omega:\wedge^2 V\to\mathbf{1}[n]$ is a morphism in $\mathcal C$. We say that it is \textit{$n$-shifted symplectic} if it is moreover non-degenerate in the following sense: $V$ is dualizable and the adjoint morphism $\omega^\flat:V\to V^*[n]$ is an equivalence\footnote{Observe that there are two candidates for being the adjoint morphism: but they may only differ by a sign as $\omega$ is skew-symmetric. The condition of being an equivalence is not affected by that ambiguity. }. 

\medskip

Note that a cochain complex $V$ is dualizable if and only if it is \textit{perfect}, meaning that its cohomology $\oplus_nH^n(V)$ is finite dimensional. Below we will consider that perfect and dualizable are synonymous. 

\begin{Exa}
When $n=0$, if $\mathcal C$ is the $\infty$-category of complexes of vector spaces, and if $V$ is concentrated in degree zero, then these notions coincide with the usual notions (pre-)symplectic vector spaces. 
\end{Exa}
\begin{Exa}[Poincar\'e duality]\label{ex1.9}
If $M$ is an $n$-dimensional oriented compact manifold, then the shifted\footnote{Without this shift, the pairing would be symmetric instead of being skew-symmetric. Indeed, recall that $S^2(V[1])=\wedge^2(V)[2]$} de Rham complex 
$DR(X)[1]=(\Omega^{*+1}(M,\mathbb{R}),d_{dR})$ is $(2-n)$-shifted symplectic when equipped with the pairing 
$\omega(\alpha,\beta):=\int_M\alpha\wedge\beta$. 
\end{Exa}

\subsubsection{Shifted symplectic structures on derived stacks}

From here we assume that $X$ is a geometric stack such that the cotangent complex $\mathbb{L}_X$ is dualizable\footnote{We will call such a derived stack an \textit{Artin stack}. }. The cotangent complex thus has a dual $\mathbb{T}_X:=\mathbb{L}_X^*$, called the \textit{tangent complex}. 

In this case a $2$-form of degree $n$ is thus a section $\omega_0$ of $S^2(\mathbb{L}_X[-1])[n+2]$, i.e.~a map $\mathcal O_X\to S^2(\mathbb{L}_X[-1])[n+2]$ in $QCoh(X)$. By duality, this is equivalent to the data of a map $S^2(\mathbb{T}_X[1])\to \mathcal O_X[n+2]$, which in turns is precisely an $n$-shifted pre-symplectic structure $\wedge^2\mathbb{T}_X\to\mathcal O_X[n]$ on $\mathbb{T}_X$ in $QCoh(X)$. 
\begin{Def}
An $n$-shifted symplectic structure on $X$ is a closed $2$-form $\omega$ of degree $n$ such that $(\mathbb{T}_X,\omega_0)$ is an $n$-shifted symplectic object. In other words, we require that $\omega_0^\flat:\mathbb{T}_X\to \mathbb{L}_X[n]$ is an equivalence. 
\end{Def}

\subsubsection{Examples of shifted symplectic structures}\label{sssec-1.2.3}

Let us first go back to Example \ref{ex1.6} and check the condition a closed $2$-form of degree $1$ on a geometric $1$-stack $X=[G_0/G_1]$ must satisfy in order to define a $1$-shifted symplectic structure. We recall several facts: 
\begin{itemize}
\item the non-degeneracy condition only concerns the leading term $\omega_0\in \Omega^2(G_1)$. 
\item quasi-coherent sheaves on $X=[G_0/G_1]$ are exactly $G_1$-equivariant sheaves on $G_0$, and the property of being an equivalence is something that one can check on the underlying sheaf on $G_0$.\footnote{In other words, the pull-back functor $QCoh(X)\to QCoh(G_0)$ is conservative. }
\item the underlying sheaf of the tangent complex $\mathbb{T}_X$ is the two-term complex 
\begin{eqnarray*}
-1 & & 0 \\
\mathcal L & \overset{\rho}{\longrightarrow} & T_{G_0}
\end{eqnarray*}
of sheaves on $G_0$, where $\mathcal L$ is the Lie algebroid of the groupoid $G_1\rightrightarrows G_0$, and $\rho$ is the anchor map\footnote{Recall that, as a vector bundle on $G_0$, $\mathcal L$ is the restriction to $G_0$ of the bundle $T^s_{G_1}$ of vectors tangent to the source map $s:G_1\to G_0$. The anchor map is given by the tangent to the target map $t:G_1\to G_0$. This calculation of the underlying sheaf $p^*\mathbb{T}_X$ of $\mathbb{T}_X$, where $p:G_0\to X=[G_0/G_1]$ is the quotient map, actually follows from a smooth descent argument, which we sketch now. Since 
$$
X=colim_{[n]\in\Delta^{op}}(\underbrace{G_1\times_{G_0}\cdots\times_{G_0}G_1}_{n\textrm{ times}})
$$
then, denoting $e:G_0\to G_1$ the identity map, 
$$
p^*\mathbb{T}_X=colim_{[n]\in\Delta^{op}}(\underbrace{e^*\mathbb{T}_{G_1}\oplus_{\mathbb{T}_{G_0}}\cdots\oplus_{\mathbb{T}_{G_0}}e^*\mathbb{T}_{G_1}}_{n\textrm{ times}})\,.
$$
Hence $p^*\mathbb{T}_X$ can be obtained as the normalized complex of the above simplicial diagram of vector bundles on $G_0$, which can be shown to be the two-term complex $\mathcal L\overset{\rho}{\longrightarrow}T_{G_0}$. }
\item $e^*T_{G_1}=\mathcal L\oplus T_{G_0}$, and $\omega_0$ is compatible with this decomposition ($e:G_0\to G_1$ is the unit map). 
\end{itemize}
Therefore $\omega_0^\flat$ induces a morphism of two-term complexes of sheaves on $X$, as follows: 
$$
\xymatrix{
\mathbb{T}_X\ar[d]&(\mathcal L \ar[r]\ar[d] & T_{G_0}\ar[d]) \\
\mathbb{L}_X[1]&(T_{G_0}^*\ar[r] & \mathcal L^*)}
$$
By duality, to check that it is a quasi-isomorphism is equivalent to check that its kernel is acyclic. 
This amounts to require that at every point $x\in G_0$, the map $\ker(\omega_{0,x})\cap\mathcal L_x\to\ker(\omega_{0,x})\cap T_xG_0$ is an isomorphism. This is precisely the non-degeneracy condition appearing in the definition of a \textit{quasi-symplectic groupoid} (see \cite[Definition 2.5]{Xu}). 

\begin{Rem}
Genuine symplectic groupoids from \cite{We2} correspond to the situation when the $3$-form $\omega_1$ (encoding the closeness of $\omega_0$ up to homotopy) is $0$. We provide a nice interpretation of this condition in terms of Lagrangian structures in the next subsection. 
\end{Rem}

\begin{Rem}
The above discussion shows in particular that, even if one deals with non-derived stacks, one has to consider cohomologically shifted $2$-forms if one wants to grasp any reasonable kind of non-degeneracy property. Indeed, the tangent complex of $[G_0/G_1]$ sits in degree $-1$ and $0$ while the cotangent complex (its dual) sits in degree $0$ and $1$. 
\end{Rem}

\begin{Exa}\label{exa-coadj}
Let $G$ be a Lie group, with Lie algebra $\mathfrak g$, and consider the groupoid $G\times\mathfrak{g}^*\rightrightarrows\mathfrak{g}^*$ of the coadjoint action (of $G$ on $\mathfrak g^*$). It is a symplectic groupoid : indeed, $G\times\mathfrak{g}^*\cong T^*G$ carries a canonical symplectic form $\omega_0$. Hence its quotient stack $[\mathfrak g^*/G]$ is $1$-shifted symplectic. There is an explicit \textit{ad hoc} description of this $1$-shifted symplectic structure (see \cite[\S1.2.3]{Cal15} and \cite[Example 2.10]{Cal14}), which we now briefly sketch. As we have seen, the tangent complex of $[\mathfrak g^*/G]$ is the two-term complex of $G$-equivariant vector bundles 
\begin{eqnarray*}
-1 & & 0 \\
\mathfrak g\otimes\mathcal O_{\mathfrak g^*} & \longrightarrow & T_{\mathfrak g^*}\cong \mathfrak g^*\otimes\mathcal O_{\mathfrak g^*}
\end{eqnarray*}
on $\mathfrak g^*$, and the cotangent complex is 
\begin{eqnarray*}
0 & & 1 \\
\mathfrak g\otimes\mathcal O_{\mathfrak g^*}\cong T^*_{\mathfrak g^*} & \longrightarrow & \mathfrak g^*\otimes\mathcal O_{\mathfrak g^*}
\end{eqnarray*}
We thus see that the tangent complex and the $1$-shifted cotangent complex are canonically identified, \textit{via} the $1$-shifted two-form 
$$
\omega=\sum_i\xi^id_{dR}x_i\in\big(\mathfrak g^*\otimes\Omega^1(\mathfrak g^*)\big)^G\,,
$$
where $(x_i)_i$ is a basis of $\mathfrak g$ (which defines coordinates on $\mathfrak g^*$) and $(\xi^i)_i$ is the dual basis of $\mathfrak g^*$. One easily see that this form is (strictly) $d_{dR}$-closed. 
\end{Exa}

\begin{Exa}
If $G$ is a Lie group equipped with an invariant non-degenerate symmetric bilinear form on its Lie algebra, then the groupoid of the conjugation action (of $G$ on itself) is quasi-symplectic, after \cite[Proposition 2.8]{Xu}. Therefore its quotient stack $[G/G]$ is $1$-shifted symplectic. We refer to \cite[\S1.2.5]{Cal15} and \cite{Saf16} for more details about this $1$-shifted symplectic structure. 
\end{Exa}

\begin{Exa}
It has also been shown in \cite{PTVV} that, if $G$ is a Lie group equipped with a $G$-invariant symmetric bilinear form $c\in S^2(\mathfrak g^*)^G$ on its Lie algebra, then the induced closed $2$-form of degree $2$ on the classifying stack $BG=[*/G]$ (see Example \ref{ex-bg}) is non-degenerate if and only if $c$ is non-degenerate (in the usual sense). We will see later that the $1$-shifted symplectic structure on $[G/G]$ can be recovered from this $2$-shifted symplectic structure on $BG$. 
\end{Exa}

\begin{Exa}
One can easily see that a symplectic $2$-groupoid $(X_\bullet,\omega)$ in the sense of \cite[Definition 2.7]{MT} 
induces a $2$-shifted symplectic structure on the geometric stack $X=colim_{[n]\in\Delta^{op}} X_n$. 
As noticed by the author of \cite{MT}, their notion of symplectic $2$-groupoid is not Morita invariant\footnote{Indeed, both their notions of closeness and non-degeneracy are too strict from the perspective of shifted symplectic structures. }.  
\end{Exa}

\subsection{Lagrangian structures: definition and examples}

\subsubsection{Lagrangian structures in the linear setting}

Let $(V,\omega)$ be a symplectic vector space, and recall that a \textit{Lagrangian} in $V$ is a subspace $L\subset V$ such that: 
\begin{itemize}
\item $L$ is \textit{isotropic}: $\omega_{|L}=0$. 
\item $L$ is non-degenerate in the sense that it is maximal. 
\end{itemize}
The isotropy condition is equivalent to the fact that the inclusion $L\subset V$ factors through $L\subset L^\circ\subset V$, where 
$L^\circ:=\{v\in V|\omega(v,L)=0\}$. 
Observe that since $\omega$ is non-degenerate\footnote{Meaning that $\omega^\flat:V\to V^*$ is an isomorphism.}, we have a canonical identification 
$$
L^\circ\tilde\longrightarrow L^\perp:=\{\xi\in V^*|\xi_{|L}=0\}
$$
with the conormal $L^\perp$ to $L$. 

The non-degeneracy condition for $L$ is equivalent to any of the following: 
\begin{itemize}
\item $dim(L)=\frac12 dim(V)$. 
\item the inclusion $L\subset L^\circ$ is an equality. 
\item the map $L\to L^\perp$ given by $\ell\mapsto\omega^\flat(\ell)$ is an isomorphism\footnote{Observe that $L^\perp$ is by definition the conormal of $L\subset V$. This equivalent characterization of non-degeneracy corresponds to the one we gave in the beginning of Section \ref{section1}. }. 
\end{itemize}
Even though this is rather unusual, the last two characterizations can be reformulated as follows: 
\begin{itemize}
\item the sequence $0\to L\to V\to L^*\to 0$, where the map $V\to L^*$ sends $v$ to $\omega^\flat(v)_{|L}$, is exact. 
\item the sequence $0\to L\to V^*\to L^*\to 0$, where the map $L\to V^*$ is given by $\ell\mapsto\omega^\flat(\ell)$, is exact. 
\end{itemize}
Note that the two sequences we have written are dual to each other: it is thus clear that one of them is exact if and only if the other is. 

This leads us to the following $\infty$-categorical generalization, which one can guess by observing that a short exact sequence of $0\to A\to B\to C\to 0$ of vector spaces is the same as a bicartesian square 
$$
\xymatrix{
A \ar[r]\ar[d] & B \ar[d]\\
0 \ar[r] & C
}
$$
\begin{Def}\label{def-lag}
Let $(\mathcal C,\otimes,\mathbf{1})$ be a stable symmetric monoidal $\infty$-category\footnote{As before, the unaccustomed reader can think about cochain complexes up to quasi-isomorphisms. }, and let $(V,\omega)$ be an $n$-shifted pre-symplectic object. \\
\textbf{a)} An isotropic structure on a morphism $L\to V$ is an homotopy between $\omega_{|L}$ and $0$. \\
\textbf{b)} An isotropic structure $\gamma$ is non-degenerate or Lagrangian if $L$ is perfect and the induced homotopy commuting square 
$$
\xymatrix{
L \ar[r]\ar[d] & V\ar[d] \\
0 \ar[r] & L^*[n]
}
$$
is (co)cartesian\footnote{First note that we are in a stable $\infty$-category, so that cartesian squares are cocartesian, and \textit{vice-versa}. The reader who is familiar with triangulated categories can think of such a square as a distinguished triangle. Readers who are not familiar neither with stable $\infty$-categories nor with triangulated categories can think of such a square as inducing a long exact sequence in cohomology $\cdots\to H^k(L)\to H^k(V)\to H^{k+n}(L^*)\to H^{k+1}(L)\to\cdots$. }
\end{Def}
\begin{Rem}
The vertical map $V\to L^*[n]$ in the above diagram is adjoint to the composed map $V\otimes L\to V\otimes V\overset{\omega}{\to} \mathbf{1}[n]$. A consequence of the definition is that $V$ is also perfect, and the non-degeneracy condition is equivalent to ask the (shifted) dual homotopy commuting square 
$$
\xymatrix{
L \ar[r]\ar[d] & V^*[n]\ar[d] \\
0 \ar[r] & L^*[n]
}
$$
to be (co)cartesian. 
Furthermore (see \cite[Lemma 1.3]{Cal17}), any $n$-shifted pre-symplectic object that has a Lagrangian is automatically $n$-shifted symplectic (i.e.~$\omega$ is non-degenerate if $\gamma$ is). 
\end{Rem}

Below we give two examples when $\mathcal C$ is the $\infty$-category of complexes of vector spaces. 

\begin{Exa}\label{ex-1.16}
Of course, ordinary Lagrangian subspaces in ordinary symplectic vector spaces are examples of Lagrangian structures for $0$-shifted symplectic. Indeed, the inclusion map carries a Lagrangian structure (which is just the constant self-homotopy of $0$). 
\end{Exa}

\begin{Exa}[Relative Poincar\'e duality]
Let $N$ be an oriented $(n+1)$-dimensional compact manifold with oriented boundary $\partial N=M$. 
Recall from Example \ref{ex1.9} that $\Omega^{*+1}(M)$ is $(2-n)$-shifted symplectic. 
We claim that the pull-back $\iota^*\Omega^{*+1}(N)\to \Omega^{*+1}(M)$ along the boundary inclusion $\iota:M\hookrightarrow N$ carries a Lagrangian structure. The isotropic structure comes from Stokes' formula: 
$$
\int_M\iota^*(\alpha\wedge\beta)=\int_Nd_{dR}(\alpha\wedge\beta)\,.
$$
The homotopy $\gamma$ is given by $\alpha\wedge\beta\mapsto \int_N\alpha\wedge\beta$. 
\end{Exa}

\subsubsection{Lagrangian structures on derived stacks}

Let $f:L\to X$ be a morphism of derived stacks and assume that $X$ carries an $n$-shifted pre-symplectic structure $\omega$. 
An \textit{isotropic} structure on $f$ (or, abusing language, ``on $L$'') is a path $\gamma$ between $f^*\omega$ and $0$ in the space $\mathcal{A}^{2,cl}(L,n)$ of closed $2$-forms of degree $n$ on $L$. 

In terms of the cocycle characterization of closed forms, this can be understdood as follows: $\gamma=\gamma_0+\gamma_1+\cdots$ with 
\begin{itemize}
\item $\gamma_i$ has weight $2+i$ and degree $1+n$. 
\item $(\partial+d_{dR})(\gamma)=f^*\omega$, meaning that $\partial\gamma_0=f^*\omega_0$, $\partial(\gamma_1)+d_{dR}(\gamma_0)=f^*\omega_1$, etc$\dots$
\end{itemize}
\begin{center}
\begin{tikzpicture}
\node[text justified] (A0) at (0,0) {$1+n$};
\node[text justified] (B0) at (0,1.5) {$2+n$};
\node[text justified] (C0) at (0,2.5) {degree$\backslash$weight};
\node[text justified] (A1) at (2,0) {$\gamma_0$};
\node[text justified] (B1) at (2,1.5) {$f^*\omega_0$};
\node[text justified] (C1) at (2,2.5) {$2$};
\node[text justified] (A2) at (4,0) {$\gamma_1$};
\node[text justified] (B2) at (4,1.5) {$f^*\omega_1$};
\node[text justified] (C2) at (4,2.5) {$3$};
\node[text justified] (A3) at (6,0) {$\gamma_2$};
\node[text justified] (B3) at (6,1.5) {$f^*\omega_2$};
\node[text justified] (C3) at (6,2.5) {$4$};
\node[text justified] (A4) at (8,0) {$\cdots$};
\node[text justified] (B4) at (8,1.5) {$\cdots$};
\node[text justified] (C4) at (8,2.5) {$\cdots$};
\draw [>=latex,->] (A1) to node[midway,fill=white]{$\partial$} (B1);
\draw [>=latex,->] (A2) to node[midway,fill=white]{$\partial$} (B2);
\draw [>=latex,->] (A3) to node[midway,fill=white]{$\partial$} (B3);
\draw [>=latex,->] (A1) to node[midway,fill=white]{$d_{dR}$} (B2);
\draw [>=latex,->] (A2) to node[midway,fill=white]{$d_{dR}$} (B3);
\draw [>=latex,->] (A3) to node[midway,fill=white]{$d_{dR}$} (B4);
\draw [>=latex,-] (-1.2,2) to (8.5,2); 
\draw [>=latex,-] (1.5,-0.2) to (1.5,2.7); 
\end{tikzpicture}
\end{center}
Let us assume that both $L$ and $X$ are Artin stacks and that $\omega$ is non-degenerate. An isotropic structure $\gamma$ of $f$ is \emph{Lagrangian} if the leading term $\gamma_0$, which can be viewed as an isotropic structure on the morphism $\mathbb{T}_L\to f^*\mathbb{T}_X$ in the sense of Definition \ref{def-lag}, is non-degenerate\footnote{Observe that one could define non-degenerate isotropic structures on $f$ even without assuming that $\omega$ is non-degenerate. These would nevertheless not deserve to be called Lagrangian structures as it would not imply that $\omega$ is symplectic (it would only imply that $f^*\omega$ is non-degenerate). }. 
\begin{Rem}
Having an isotropic structure on $\mathbb{T}_L\to f^*\mathbb{T}_X$ tells us that we have a morphism from $\mathbb{T}_L$ to the homotopy fiber of $f^*\mathbb{L}_X[n]\to \mathbb{L}_L[n]$, which is nothing but $\mathbb{L}_f[n-1]$. The non-degeneracy condition tells us then that this morphism $\mathbb{T}_L\to\mathbb{L}_f[n-1]$ is an equivalence. In the case when $n=0$ and the map $L\to X$ is an inclusion of genuine manifolds, then $\mathbb{T}_L$ is the usual tangent bundle $TL$ and $\mathbb{L}_f[-1]$ is the conormal bundle $N^*L$. Hence it coincides with the usual notion of a Lagrangian submanifold that we recalled at the beginning of Section \ref{section1}. 
\end{Rem}
Observe that the notion of a Lagrangian structure on a morphism fully addresses problem \textbf{[d]} from the Introduction. 

\subsubsection{Examples of Lagrangian structures}

\begin{Exa}
Of course, any genuine Lagrangian submanifold in a genuine symplectic manifold provide an example of a Lagrangian structure, on the inclusion morphism. Indeed, as in Example \ref{ex-1.16} the Lagrangian structure is the constant self-homotopy of the closed $2$-form $0$. 
\end{Exa}

\begin{Exa}[Symplectic is Lagrangian]\label{ex-symplag}
Let $*_{(n)}$ be the point equipped with the trivial $n$-shifted symplectic structure, given by the zero $2$-form. Surprisingly enough, it was noticed in \cite[Example 2.3]{Cal15} that the space of Lagrangian structures on the morphism $X\to *$ is equivalent to the space of $(n-1)$-shifted symplectic structures on $X$. This answers part of question \textbf{[b]} in the Introduction: a Lagrangian ``in'' $P_{red}=*$ is a $(-1)$-shifted symplectic stack.  
\end{Exa}

\begin{Exa}\label{ex-conormor}
If $X\to Y$ is a morphism of Artin stacks, then it has been shown in \cite{Cal17} that the cotangent stack $T^*Y$ is $0$-shifted symplectic and that the morphism $T^*_XY\to T^*Y$ from the conoraml stack to the cotangent stack is Lagrangian. 
\end{Exa}

\begin{Exa}[Symplectic groupoids]\label{ex-1.20}
Let $G_1\rightrightarrows G_0$ be a Lie groupoid together with a quasi-symplectic structure $\omega=\omega_0+\omega_1$ in the sense of \cite[Definition 2.5]{Xu}. We have seen above in paragraph \ref{sssec-1.2.3} that the stack $[G_0/G_1]$ then carries a $1$-shifted symplectic structure. 
Notice that the pull-back of $\omega$ along the quotient map $G_0\to [G_0/G_1]$ is $\omega_1$. Hence the quotient map carries a Lagrangian structure if and only if $\omega_1=0$, meaning that $G_1\rightrightarrows G_0$ actually is a symplectic groupoid. 
\end{Exa}
\begin{Rem}
One actually has an equivalence between the following three sets of data: 
\begin{itemize}
\item a $1$-shifted symplectic structure on the quotient $[G_0/G_1]$ and a Lagrangian structure on the quotient map $G_0\to [G_0/G_1]$. 
\item a quasi-symplectic structure on $G_1\rightrightarrows G_0$ is such that $\omega_1=0$ (i.e.~a multiplicative symplectic structure on $G_1$). 
\item a symplectic structure on $G_1$ such that the submanifold $G_0\subset G_1$ of units of the groupoid is Lagrangian. 
\end{itemize}
This can been seen in the case of the coadjoint action groupoid $G\times\mathfrak g^*\rightrightarrows \mathfrak g^*$ from Example \ref{exa-coadj}:  
\begin{itemize}
\item pulling back the $1$-shifted two-form $\sum_i\xi^id_{dR}x_i$ along $\mathfrak g^*\to [\mathfrak g^*/G]$ amounts to setting $\xi^i=0$, and thus gives $0$. 
\item the symplectic structure on $G\times\mathfrak g^*\cong T^*G$ is multiplicative. 
\item $\{e\}\times\mathfrak g^*\subset G\times\mathfrak g^*\cong T^*G$ is Lagrangian. 
\end{itemize}
\end{Rem}

\begin{Exa}[Hamiltonian groupoid actions]\label{ex-xuxu}
Let $(G_1\rightrightarrows G_0,\omega=\omega_0+\omega_1)$ be a quasi-symplectic groupoid and let $J:X_0\to G_0$ be a $G_1$-space, that is to say a family over $G_0$ together with a $G_1$-action. 
There is a nice description of Lagrangian structures on the map $[J]:[X_0/G_1]\to [G_0/G_1]$. 
First of all, let $X_\bullet$ be the nerve of the action groupoid $X_1\rightrightarrows X_0$, where $X_1=G_1\times_{G_0}X_0$. 
Then the map $[J]$ is presented by the obvious morphism of simplicial manifolds $J_\bullet:X_\bullet\to G_\bullet$ (note that $J_0=J$). 
Then the pull-back $[J]^*\omega$ is $J^*\omega_1+J_1^*\omega_0$. The condition that $[J]$ is isotropic therefore reads as follows: there exists a $2$-form $\eta\in\Omega^2(X_0)$ such that $[J]^*\omega=\partial\eta+d_{dR}\eta$. In other words: 
$$
J^*\omega_1=d_{dR}\eta\qquad\mathrm{and}\qquad J_1^*\omega_0=\partial\eta\,.
$$
The second condition says that the graph of the action in $G_1\times X_0\times X_0$ is isotropic with respect to the $2$-form 
$(\omega_0,\eta,-\eta)$. Hence $[J]$ is isotropic if and only if $X_0$ is a pre-Hamiltonian $G_1$-space in the sense of \cite[Definition 3.1]{Xu}. It is an exercise to check that the isotropic structure given by $\eta$ is Lagrangian if and only if it turns $X_0$ into a Hamiltonian $G_1$-space in the sense of \cite[Definition 3.5]{Xu}. 
\end{Exa}

\begin{Exa}[Moment maps]\label{ex-moment}
Hamiltonian spaces for the symplectic groupoid $G\times\mathfrak g^*\rightrightarrows\mathfrak g^*$ are symplectic $G$-spaces $X$ together with a map $\mu:X\to \mathfrak g^*$ satisfying the moment condition $\omega^\flat(\vec{x})=d\mu_x$ for every $x\in\mathfrak{g}$. According to the previous Example, we recover the fact that the induced map $[X/G]\to [\mathfrak g^*/G]$ carries a Lagrangian structure (see e.g.~\cite[\S2.2.1]{Cal15} or \cite[Example 2.14]{Cal14}). Indeed, the moment condition implies that 
$$
[\mu]^*\sum_i\xi^id_{dR}x_i=\sum_i\xi^i\mu^*d_{dR}=\sum_i\xi^id\mu_{x_i}=\sum_i\xi^id\omega_X^\flat(\vec{x_i})=\partial\omega_X
$$
\end{Exa}

\begin{Exa}[Lie group valued moment maps]
Hamiltonian spaces for the conjugation quasi-symplectic groupoid $G\times G\rightrightarrows G$ are precisely the quasi-Hamiltonian $G$-spaces of Alekseev-Malkin-Meinrenken \cite{AMM}, the map $J:X\to G$ being called a Lie group valued moment map. 
It again follows from Example \ref{ex-xuxu} that the induced map $[X/G]\to [G/G]$ then carries a Lagrangian structure (see e.g.~\cite[\S2.2.2]{Cal15} or \cite[\S2.3]{Saf16}). 
\end{Exa}

\begin{Exa}[Symplectic groupoids again]
If we are given a symplectic groupoid $(G_1\rightrightarrows G_0,\omega_0)$ then $G_1$ naturally becomes a Hamiltonian $G_1$-space for the action by left multiplication. In particular, we recover the fact that, for symplectic groupoids, the map $G_0=[G_1/G_1]\to [G_0/G_1]$ carries a Lagrangian structure (see Example \ref{ex-1.20}).  
\end{Exa}

\section{Derived interpretation of classical constructions}

\subsection{Weinstein's symplectic category}\label{weinstein-cat}

The cotangent space $X=T^*M$ of a differentiable manifold $M$ is symplectic. 
Let us look at it as the phase space of a classical mechanical system: $T^*M$ the space of all possible pairs of position and momentum, and we have a Hamiltonian function $H:X\to\mathbb{R}$ leading to a Hamiltonian vector field $\mathfrak{X}_H$ on $X$. 
Inside the \textit{space of histories}, here the space of paths $\gamma:[0,1]\to X$, 
one is interested in the solution space $\mathcal S$ of the equations of motions
$$
\dot{\gamma}(t)=\mathfrak{X}_{H,\gamma(t)}\,.
$$
It happens to be isomorphic to $X$ itself, as a solution is uniquely determined by the initial condition $\gamma(0)\in X$. 
The evaluation map sending $\gamma$ to the pair $\big(\gamma(0),\gamma(1)\big)$ exhibits $\mathcal S$ as a Lagrangian 
submanifold of $X\times \overline{X}$, where $\overline{X}$ stands for $X$ equipped with the opposite symplectic structure. 

This is actually a special case of the graph of a symplectomorphism\footnote{In the above example, the symplectomorphism is the time $1$ flow of the Hamiltonian vector field $\mathfrak{X}_H$. } ${\varphi:X\to Y}$ being a Lagrangian 
submanifold in $X\times\overline{Y}$. The graph of the composition of two symplectomorphisms actually coincides with the composition of their associated correspondences. Elaborating on this observation, Weinstein suggested in \cite{We} to construct a category with objects being symplectic manifolds and morphisms being Lagrangian correspondences\footnote{This idea is at the origin of the definition of a symplectic groupoid, where the graph of the multiplication is required to be Lagrangian. }. The main problem being that this is not quite a category, since there are well-known transversality issues when it comes to compose Lagrangian correspondences. These issues are resolved in the derived framework. 

Roughly speaking, there is an $\infty$-category $\mathcal{L}ag_n$ where objects are $n$-shifted symplectic stacks ($n$ being fixed) and morphisms being Lagrangian correspondences $L\to X\times \overline{Y}$. One can compose Lagrangian correspondences (see \cite[Theorem 4.4]{Cal15}): given Lagrangian correspondences $L_1\to X\times \overline{Y}$ and $L_2\to Y\times \overline{Z}$, the (derived) fiber product $L_1\times_YL_2\to X\times \overline{Z}$ is again a Lagrangian correspondence: 
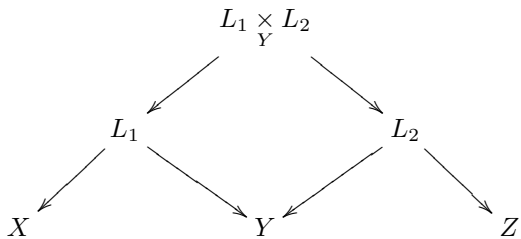
\[
  \xymatrix{
    && L_1\underset{Y}{\times}L_2 \ar[dl] \ar[dr]\\
    & L_1 \ar[dl] \ar[dr] && L_2 \ar[dl] \ar[dr] \\
    X && Y  && Z
  }
\]
The full construction of this $\infty$-category has been achieved in \cite[Section 11]{Hau}. 

\begin{Rem}
Applying the above to the case when $X=Z=*$, and combining it with Example \ref{ex-symplag}, one recovers the striking observation from \cite{PTVV} that the (derived) fiber product $L_1\times_Y L_2$ of two Lagrangian morphisms $L_i\to Y$ is $(n-1)$-shifted symplectic (if we started with an $n$-shifted symplectic $Y$).  
\end{Rem}

\subsubsection{Examples}

Below we give several examples of Lagrangian correspondences and compositions of them. 

\begin{Exa}[Lagrangian morphisms]
Every Lagrangian morphism $L\to X$ can be viewed as a Lagrangian correspondence between $X$ and $*$ (in whichever direction): 
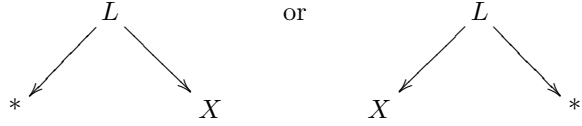
\[
  \xymatrix{
    & L \ar[dl] \ar[dr] & \\
    \textrm{*} & & X
  }
\qquad\mathrm{or}\qquad 
  \xymatrix{
    & L \ar[dl] \ar[dr] & \\
    X & & \textrm{*}
  }
\]
\end{Exa}

\begin{Exa}[Conormal to a graph]
The conormal $N^*f\cong X\times_YT^*Y\to T^*X\times\overline{T^*Y}$ to the graph $X\to X\times Y$ of a morphism $f:X\to Y$ naturally carries a shifted Lagrangian structure. This construction is actually functorial: if $g:Y\to Z$ is another morphism we then get that 
$$
N^*(g\circ f) \cong N^*f\times_{T^*Y}N^*g
$$
as Lagrangian correspondences from $T^*X$ to $T^*Z$. Indeed, the composition of correspondences 
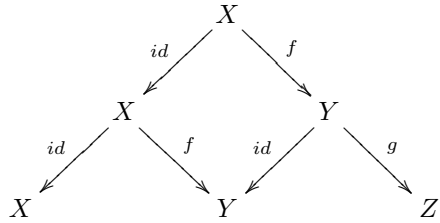
\[
  \xymatrix{
    && X \ar[dl]_{id} \ar[dr]^{f}\\
    & X \ar[dl]_{id} \ar[dr]^{f} && Y \ar[dl]_{id} \ar[dr]^{g} \\
    X && Y  && Z
  }
\]
is sent by $N^*$ to the composition of Lagrangian correspondences 
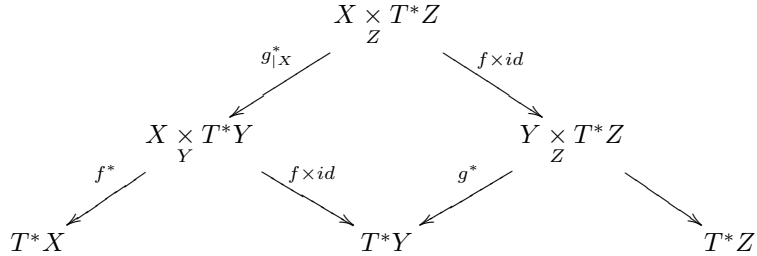
\[
\xymatrix{
    && X\underset{Z}{\times}T^*Z \ar[dl]_{g^*_{|X}} \ar[dr]^{f\times id} \\
    & X\underset{Y}{\times}T^*Y \ar[dl]_{f^*} \ar[dr]^{f\times id} && Y\underset{Z}{\times}T^*Z \ar[dl]_{g^*} \ar[dr] \\
    T^*X && T^*Y  && T^*Z
  }
\]
We thus have a functor $N^*:d\mathcal{A}rt\mathcal{S}t\to \mathcal Lag_0$ of $\infty$-categories from derived Artin stacks (with their usual morphisms) to $0$-shifted symplectic stacks with Lagrangian correspondences. 
\end{Exa}

\begin{Exa}[Lagrange multipliers/Constrained critical locus]\label{ex-lagmul}
Following \cite{RS}, we consider a \textit{variational family}, being the data of a morphism\footnote{As we work in the derived setting, and contrary to \cite{RS}, we don't require that the morphism is a surjective submersion. } $f:P\to X$ and a function $S:P\to\mathbb{R}$. 

First of all observe that the $1$-form $dS:P\to T^*P$ carries a Lagrangian structure. This is for instance proven in \cite[\S2.4]{Cal17} (this is actually true for every closed $1$-form), but it can also be obtained as a composition of Lagrangian correspondences: indeed, it is the composition of $N^*S\to T^*P\times\overline{T^*\mathbb{R}}$ with the Lagrangian $\mathbb{R}\hookrightarrow T^*\mathbb{R}\cong\mathbb{R}^2$ given by $x\mapsto(x,1)$: 
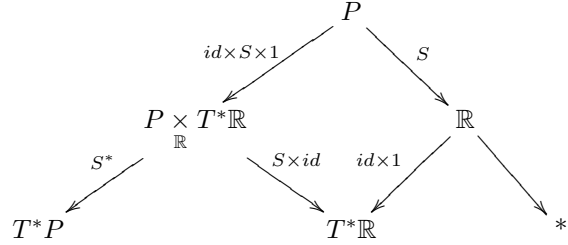
\[
\xymatrix{
    && P \ar[dl]_{id\times S\times 1} \ar[dr]^{S}\\
    & P\underset{\mathbb{R}}{\times}T^*\mathbb{R} \ar[dl]_{S^*} \ar[dr]^{S\times id} && \mathbb{R} \ar[dl]_{id\times 1} \ar[dr] \\
    T^*P && T^*\mathbb{R}  && \textrm{*}
  }
\]
Then consider the \textit{derived constrained critical locus} (or \textit{derived fiber critical locus}, or \textit{derived Lagrange multiplier space}) $\mathbf{Crit}_f(S):=P\times_{T^*P}N^*f$ of the family $f$. The morphism $\mathbf{Crit}_f(S)\to T^*X$ therefore carries a natural Lagrangian structure\footnote{Somehow generalizing \cite[Proposition 1.1]{RS}. }: 
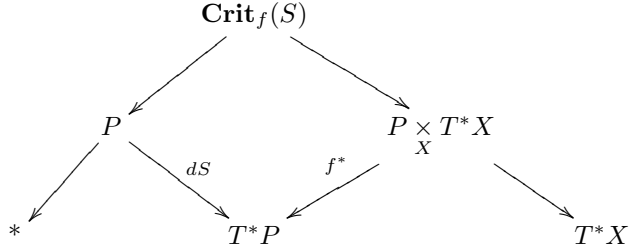
\[
\xymatrix{
    && \mathbf{Crit}_f(S) \ar[dl] \ar[dr] \\
    & P \ar[dl] \ar[dr]^{dS} && P\underset{X}{\times}T^*X  \ar[dl]_{f^*} \ar[dr] \\
    \textrm{*} && T^*P  && T^*X
  }
\]
\end{Exa}

\begin{Exa}[Derived critical locus]
Notice that when $X=*$ in the above Example then $\mathbf{Crit}_f(S)=\mathbf{Crit}(S)$ is the (absolute) derived critical locus of $f$ and is thus $(-1)$-shifted symplectic (according to Example \ref{ex-symplag}). 
\end{Exa}

\begin{Exa}[Hamiltonian bimodules]\label{ex-bimod}
Clearly, given two quasi-symplectic groupoids $G_1\rightrightarrows G_0$ and $H_1\rightrightarrows H_0$, and a Hamiltonian bimodule $X$ (in the sense of \cite[Definition 3.13]{Xu}) for these, then $[X/G_1\times H_1^{op}]$ provides a Lagrangian correspondence between $[G_0/G_1]$ and $[H_0/H_1]$. When the composition of two such bimodules is well-defined (such as in \cite[Theorem 3.16]{Xu} for instance) then the resulting Lagrangian correspondence is the composition of the Lagrangian correspondences associated with these two Hamiltonian bimodules. 

A particular case of these compositions is when we have two Hamiltonian $G_1$-spaces $X_0$ and $Y_0$: the derived fiber product 
$[X_0/G_1]\times_{[G_0/G_1]}[Y_0/G_1]\cong\big[X_0\times_{G_0} Y_0/G_1\big]$ is then $0$-shifted symplectic (again, this has to be compared with \cite[Theorem 3.21]{Xu}). Recalling from Example \ref{ex-1.20} that we also have a Lagrangian morphism $G_0\to[G_0/G_1]$, so that we get the following commuting cube of cartesian squares: 
\[
 \xymatrix {
    X_0\times_{G_0}Y_0 \ar[rr] \ar[dd] \ar[dr] && Y_0 \ar[dr] \ar[dd] |!{[dl];[dr]}\hole \\
    &\big[X_0\times_{G_0}Y_0/G_1\big] \ar[rr] \ar[dd] && [Y_0/G_1] \ar[dd]^{Lag} \\
    X_0 \ar[rr] |!{[ur];[dr]}\hole \ar[dr] && G_0 \ar[rd]_{Lag} \\
    & [X_0/G_1] \ar[rr]_{Lag} && [G_0/G_1] \\
  }
\]
\begin{itemize}
\item There are three Lagrangian morphisms to the $1$-shifted symplectic stack $[G_0/G_1]$. 
\item Their three derived intersections give rise to three $0$-shifted symplectic stacks. 
\item $X_0\times_{G_0}Y_0$ realizes three Lagrangian correspondences (for each pair of these $0$-shifted symplectic stacks). 
\end{itemize}
\end{Exa}


\subsubsection{Hamiltonian reduction}

Let us consider a phase space (i.e.~a symplectic manifold) $X$ having a Lie group of symmetries $G$ and moments for these symmetries: for each infinitesimal generator $x\in \mathfrak g$ of the action there is a Hamiltonian function $\mu_x$, i.e. $\omega^\flat(\vec{x})=d\mu_x$. Collecting all $\mu_x$'s we get a map $\mu:X\to g^*$, called a \textit{moment map}. 
 
For every weakly regular value $\xi\in\mathfrak{g}^*$ satisfying nice enough hypotheses (see \cite[Theorem 1]{MW}), the so-called \textit{reduced space} $\mu^{-1}(\mathcal O_\xi)/G$ is symplectic, where $\mathcal O_\xi$ is the coadjoint orbit of $\xi$. 
Observe that these hypotheses actually guarantee that 
the reduced space $\mu^{-1}(\mathcal O_\xi)/G$ is equivalent to the \textit{derived reduced space} 
$$
\big[X\times_{\mathfrak g^*}\mathcal O_\xi/G\big]\cong[X/G]\times_{[\mathfrak g^*/G]}[\mathcal O_\xi/G]\,.
$$

Now recall from Example \ref{ex-moment} that the map $[\mu]:[X/G]\to [\mathfrak g^*/G]$ carries a Lagrangian structure and that the map $[\mathcal O_\mu/G]\to [\mathfrak g^*/G]$ does as well\footnote{Indeed, the inclusion of coadjoint orbit $\mathcal O_\mu\hookrightarrow \mathfrak g^*$ is a moment map. }. Hence we get that the derived reduced space is naturally $0$-shifted symplectic (being a ``Lagrangian intersection in'' a $1$-shifted symplectic stack) . 

\medskip

As the reader may already have noticed, this is a particular instance of the very general phenomenon that we have already pointed in Example \ref{ex-bimod}. We thus again have a commuting cube of cartesian squares as in Example \ref{ex-bimod}: 

\[
 \xymatrix {
    X\times_{\mathfrak g^*}\mathcal O_\xi \ar[rr] \ar[dd] \ar[dr] && \mathcal O_\xi \ar[dr] \ar[dd] |!{[dl];[dr]}\hole \\
    &\big[X\times_{\mathfrak g^*}\mathcal O_\xi/G\big] \ar[rr] \ar[dd] && [\mathcal O_\xi/G] \ar[dd]^{Lag} \\
    X \ar[rr] |!{[ur];[dr]}\hole \ar[dr] && \mathfrak{g}^* \ar[rd]_{Lag} \\
    & [X/G] \ar[rr]_{Lag} && [\mathfrak g^*/G] \\
  }
\]

Let us provide below a few more examples of this phenomenon. 

\begin{Exa}[Hamiltonian spaces for symplectic groupoids are symplectic]
Assume we are given a symplectic groupoid $G_1\rightrightarrows G_0$ and a Hamiltonian $G_1$-space $X$. Recall that $G_1$ itself is a Hamiltonian $G_1$-space, so that 
$$
[X/G_1]\times_{[G_0/G_1]}[G_1/G_1]=[X/G_1]\times_{[G_0/G_1]}G_0\cong X
$$
is thus $0$-shifted symplectic. The morphism $X\to G_0$ plays the role of a moment map. 
\end{Exa}

\begin{Exa}[Quasi-Hamiltonian reduction]
Consider the conjugation quasi-symplectic groupoid $G\times G\rightrightarrows G$ and a quasi-Hamiltonian $G$-space $X$. 
Then we get that for any conjugacy class\footnote{Conjugacy classes are examples of quasi-Hamiltonian $G$-spaces (see \cite[Proposition 3.1]{AMM}). } $C\subset G$ the derived fiber product $\big[X\times_GC/G\big]$ is $0$-shifted symplectic. This gives back the fact that, under suitable assumptions, we have a genuine symplectic manifold $\mu^{-1}(C)/G$ as in \cite{AMM}, where $\mu:X\to G$ is the Lie group valued moment map. 
\end{Exa}

\subsection{Transgression}

A transgression is a map that transfers cohomology classes in a way that changes the cohomological degree, typically when integrating along fibers in differential geometry. This has for instance many manifestations in field theory, where the expression of the action functional in terms of the Lagrangian can sometimes be understood as a transgression procedure\footnote{We refer to \cite{cohesive} for a systematic approach, and to the survey papers \cite{FSS,Schreiber} for shorter and less abstract expositions. Observe that these references also deal with variants of shifted pre-symplectic structures on stacks, but the non-degeneracy condition is almost never satisfied as everything takes place in the realm of \textit{underived} stacks.}.

The so-called \textit{AKSZ formalism} from \cite{AKSZ}, which allows to present many classical gauge theories as $\sigma$-models and make them fit into the BV formalism, is also entirely based on a transgression procedure. In what follows we explain a far reaching generalization of it, called \textit{PTVV formalism}\footnote{The main advantage of the PTVV formalism, compared to the AKSZ one, is that is is model independent. For instance, all notions (e.g.~closed forms) and properties (e.g.~non-degeneracy) are invariant under appropriate equivalences (e.g.~quasi-isomorphisms). }, developed in \cite{PTVV}. 

\subsubsection{Integration theory on stacks: shifted orientations}\label{integration}

In \cite{PTVV} the authors introduce a class of derived stacks that carry a good integration theory of cohomological degree $d$. 
Such a derived stack $\Sigma$ shall be: 
\begin{itemize}
\item[\textbf{(o1)}] such that, for any other stack $F$, one can extract the $(0,*)$-part of a closed form on the product $\Sigma\times F$. One thus has a morphism $DR(\Sigma\times F)\to C^\infty(\Sigma)\otimes DR(F)$ of graded mixed complexes, where $C^\infty(\Sigma):=\Gamma(\mathcal O_\Sigma)$ is the cochain complex of derived global functions on $\Sigma$. 
\item[\textbf{(o2)}] such that for any stack $F$ the derived global section functor 
$$
\Gamma(\Sigma,-):QCoh(\Sigma\times F)\to QCoh(F)
$$
preserves perfect objects\footnote{Recall that perfect objects are roughly locally equivalent to bounded complexes of finite dimensional vector bundles. This condition can thus be understood as saying that the cohomology of finite dimensional vector bundles on $Sigma$ is finite dimensional. }. 
\item[\textbf{(o3)}] equipped with a \textit{$d$-class} $[\Sigma]:C^\infty(\Sigma)\to\mathbb{R}[-d]$, also called \textit{$d$-shifted pseudo-orientation}. 
\item[\textbf{(o4)}] such that for any perfect complex $E$ on $\Sigma$, the pairing 
$$
\Gamma(\Sigma,E)\otimes\Gamma(\Sigma,E^*)[d]\to C^\infty(\Sigma)[d]\to \mathbb{R}
$$
is non-degenerate in cohomology. 
\end{itemize}
A pair $(\Sigma,[\Sigma])$ as above is called a \textit{$d$-oriented stack}. If only \textbf{(o1)} and \textbf{(o3)} are satisfied we call it a \textit{$d$-pseudo-oriented} stack. 

\begin{Exa}[The de Rham stack $M_{DR}$]\label{ex-DR}
Let $M$ be a closed compact oriented $d$-manifold and let $\Sigma=M_{DR}$ be the de Rham stack of $M$. It is the stack obtained as the quotient of the groupoid $\widehat{M\times M}\rightrightarrows M$, where $\widehat{M\times M}$ is the formal neighborhood of the diagonal in $M$ (i.e.~we identify any two infinitesimally close points). A (complex of) vector bundle on $M_{DR}$ is a (complex of) vector bundle on $M$ together with a flat connection, and one has that $C^\infty(M_{DR})\cong(\Omega^*(M),d_{dR})$ is the de Rham complex of $M$. We thus have a $d$-shifted orientation given by $\int_M:(\Omega^*(M),d_{dR})\to\mathbb{R}[-d]$. 

The de Rham stack has first been introduced by Carlos Simpson in the algebro-geometric context \cite{simpson}. 
\end{Exa}
\begin{Exa}[The Betti stack $M_B$]\label{ex-Betti}
Let $M$ be a closed compact oriented $d$-manifold and let $\Sigma=M_B$ be the Betti stack of $M$. It is the constant stack associated with the homotopy type of $M$. Given a cellular decomposition of $M$, $M_B$ can be described as the stack with $0$-cells as points, $1$-cells as isomorphisms, $\dots$, $n$-cells as $n$-isomorphisms, etc$\dots$ One has that $C^\infty(M_B)\cong C^*_{sing}(M,\mathbb{R})$ is the singular cochain complex of $M$. We have a $d$-shifted orientation given by evaluating cochains on the fundamental class of $M$. 
\end{Exa}
\begin{Exa}[of a pseudo-orientation]\label{ex-curvature}
Let $M$ be a closed compact oriented surface and let $\Sigma=B(\wedge^2T_M)$ be the split first order degree $-1$ extension\footnote{For a bundle $E\to M$, the split first order degree $-1$ extension is also the classifying stack of $E$, viewed as a bundle of groups on $M$ for the $+$ law. A more standard notation in differential (super)geometry is $E[1]$. } of $M$ by the sheaf $\wedge^2T_M$: $C^\infty(\Sigma)\cong C^\infty(M)\oplus\Omega^2(M)[-1]$. We have a $1$-shifted pseudo-orientation\footnote{One could imagine a Banach framework in which this example is a $1$-shifted orientation. } given by $\int_M$. 
\end{Exa}
\begin{Rem}
Observe that $\wedge^2T_M$ gets a nice geometric interpretation as the subspace of $M^{\Delta^2_{inf}}$ consisting of those infinitesimal $2$-simplicies in $M$ that are either non-degenerate or constant. In other words, $\wedge^2TM$ is the space of ``order one surfaces'' in $M$ (like $TM$ is the space of order one curves in $M$). 
\end{Rem}

Given a $d$-pseudo-oriented derived stack $(\Sigma,[\Sigma])$ and any derived stack $F$ one can compose the morphism of graded mixed complexes $DR(\Sigma\times F)\to C^\infty(\Sigma)\otimes DR(F)$ from \textbf{(o1)} with the $d$-shifted pseudo-orientation $[\Sigma]$, and get a morphism of graded mixed complexes $DR(\Sigma\times F)\to DR(F)[-d]$, denoted $\int_{[\Sigma]}$. In particular, $\int_{[\Sigma]}$ induces a map $\mathcal A^{p,cl}(\Sigma\times F,n)\to \mathcal A^{p,cl}(F,n-d)$. 

Therefore, if one let $F$ be a derived mapping stack $F:=\mathbf{Map}(\Sigma,X)$, then we get a map $\mathcal A^{p,cl}(X,n)\to \mathcal A^{p,cl}(F,n-d)$: it is given by $\omega\mapsto\int_{[\Sigma]}ev^*\omega$, where $ev:\Sigma\times\mathbf{Map}(\Sigma,X)\to X$ is the evaluation map. 
\[
\xymatrix{
\Sigma\times \mathbf{Map}(\Sigma,X) \ar[r]^-{ev} \ar[d]& X\\
\mathbf{Map}(\Sigma,X) &
}
\qquad
\xymatrix{
\mathcal A^{p,cl}\big(\Sigma\times \mathbf{Map}(\Sigma,X)\big) \ar[d]_{\int_{[\Sigma]}} & \ar[l]_-{ev^*} \mathcal A^{p,cl}(X,n)\\
\mathcal A^{p,cl}\big(\mathbf{Map}(\Sigma,X),n-d\big) &
}
\]
\begin{Exa}
Let $\Sigma=S^1_B\cong B\mathbb{Z}$ and let $X=BG$ be the classifying stack of a Lie group $G$. Then $\mathbf{Map}(\Sigma,X)\cong[G/G]$. 
Consider the $2$-shifted pre-symplectic structure $\omega$ on $BG$ determined by an invariant symmetric pairing $c\in S^2(\mathfrak g^*)^G$. One can show (see \cite{Saf16}) that the transgressed $1$-shifted pre-symplectic form on $[G/G]$ we get is indeed the one coming from the 
pre-symplectic groupoid $G\times G\rightrightarrows G$. 
\end{Exa}

\begin{Exa}
Let $\Sigma=S^1_{DR}$ and let $X=BG$ be the classifying stack of a connected Lie group $G$. Then $\mathbf{Map}(\Sigma,X)\cong\big[\Omega^1(S^1,\mathfrak g)/C^\infty(S^1,G)\big]$, where the action is given by $g\cdot\alpha=Ad_g\alpha+dgg^{-1}$. 
One can show (by an explicit calculation) that the transgressed $1$-shifted pre-symplectic form on $[G/G]$ we get is indeed the one pulled-back from $\big[C^\infty(S^1,\mathfrak g)^*/C^\infty(S^1,G)\big]$ along the ``pairing+integration map'' $\Omega^1(S^1,\mathfrak g)\to C^\infty(S^1,\mathfrak g)^*$. 
\end{Exa}

\begin{Exa}
Let $\Sigma=M[\wedge^2T_M]$ be as in Example \ref{ex-curvature} and let $X=BG$ be the classifying stack of a compact and simply connected group $G$. Consider the $2$-shifted pre-symplectic structure $\omega$ on $BG$ determined by an invariant symmetric pairing $c\in S^2(\mathfrak g^*)^G$. Then $\int_{[\Sigma]}ev^*\omega$ is a $1$-shifted pre-symplectic structure on 
$\mathbf{Map}(\Sigma,X)\cong\big[\Omega^2(M,\mathfrak g)/C^\infty(M,G)\big]$.\footnote{For a more general $G$, $\mathbf{Map}(\Sigma,X)$ is the derived moduli stack of $G$-bundles equipped with a basic $2$-form. } 
This $1$-shifted pre-symplectic stack is an ``infinite dimensional analog'' of $[\mathfrak g^*/G]$ and will play a crucial role in the derived interpretation of the infinite dimensional reduction procedure that we have seen in the Introduction. 
\end{Exa}

We now come to a very useful result from \cite{PTVV}: 

\begin{Thm}\label{thm-ptvv}
If $(\Sigma,[\Sigma])$ is a $d$-oriented stack, $(X,\omega)$ is an $n$-shifted symplectic stack, and $\mathbf{Map}(\Sigma,X)$ is an Artin stack, then $\int_{[\Sigma]}ev^*\omega$ is non-degenerate. 
Therefore $\mathbf{Map}(\Sigma,X)$ naturally becomes $(n-d)$-shifted symplectic. 
\end{Thm}

This is a very nice statement as 
\begin{itemize}
\item many symplectic structures on moduli spaces can be recovered as particular instances of this result (for $n=d$). 
\item several examples of so-called perfect obstruction theories (after Behrend--Fantechi \cite{BF}) as well (for $n=d-1$). 
\end{itemize}

\subsubsection{De Rham stack \textit{versus} Betti stack}

Let $X$ be a $d$-dimensional manifold. There is a map $X_{DR}\to X_B$ from the de Rham stack to the Betti 
stack\footnote{For instance, consider a cover $\pi:\mathfrak U=\coprod_i U_i$ by contractible open subsets of $X$ such that all iterated intersections of these open subsets are contractible as well. 
In other words, the nerve of the map $\pi:\mathfrak U\to X$ is made of contractible open subsets of $X$. 
Note that $X$ is equivalent to the homotopy colimit (in stacks) of the nerve $N(\pi)$ of the map $\pi:\mathfrak U\to X$. 
Similarly $X_{DR}$ is equivalent to the homotopy colimit of $N(\pi)_{DR}$ (where we have applied the de Rham functor levelwise). 
Finally, using the fact that we always have a terminal map $(\mathbb{R}^d)_{DR}\to*$ we get a morphism $N(\pi)_{DR}\to \pi_0N(\pi)=X_B$ (here again we have applied $\pi_0$ levelwise).}. 
This tells us in particular that for a derived stack $F$, we have morphism $\mathbf{Map}(X_B,F)\to \mathbf{Map}(X_{DR},F)$. 
It can be shown along the lines of \cite{Porta} that, whenever $F=BG$, this map is an equivalence\footnote{In \cite{Porta} 
this is proven for $F$ being the derived stack of perfect complexes and in the context of derived analytic stacks. The proof carries over in the differentiable context, and for $F=BG$ as well. }. 
\begin{Exa}
Let $F=BG$ and $X=S^1$. On the one hand recall that $S^1_B\cong B\mathbb{Z}$ and thus $\mathbf{Map}(X_B,F)\cong [G/G]$ (where the action of $G$ on itself is by conjugacy, as usual). On the other hand one can show that $\mathbf{Map}(X_{DR},F)$ can be described as the quotient stack\footnote{Assuming $G$ is connected. } $[\Omega^1(S^1,\mathfrak g)/C^\infty(S^1,G)]$, where the action is given by $g\cdot\alpha=Ad_g\alpha+dgg^{-1}$. 
We therefore get that there is an equivalence between $[\Omega^1(S^1,\mathfrak g)/C^\infty(S^1,G)]$ and $[G/G]$. 

For a more general $X$ this tells us that the derived stack of $G$-local systems on $X$ is equivalent to the derived stack of flat $G$-connections on $X$. 
\end{Exa}

Additionally, one can prove that if $X$ is compact and oriented then the $d$-orientations on $X_{DR}$ and $X_B$ do coincide. 
This in particular tells us that, if $F$ is an $n$-shifted symplectic derived stack then we have an equivalence \textbf{of $(n-d)$-shifted symplectic stacks} $\mathbf{Map}(X_{DR},F)\cong \mathbf{Map}(X_B,F)$. 

\begin{Exa}
Going back to the previous example, and assuming now that $G$ is compact and connected, we get an equivalence of $1$-shifted symplectic stacks between $[\Omega^1(S^1,\mathfrak g)/C^\infty(S^1,G)]$ and $[G/G]$. We therefore get a correspondence between Lagrangian morphisms to $[\Omega^1(S^1,\mathfrak g)/C^\infty(S^1,G)]$ and Lagrangian morphisms to $[G/G]$, providing a very nice interpretation the correspondence between Hamiltonian $L(G)$-spaces and quasi-Hamiltonian $G$-spaces from \cite{AMM}. 

For a more general $d$-dimensional compact oriented $X$, we get that the derived stacks of $G$-local systems and of flat $G$-connections on $X$ are equivalent as $(2-d)$-shifted symplectic stacks. 
\end{Exa}

We would now like to calculate the two $(2-d)$-shifted symplectic structures from the above example at an $\mathbb{R}$-point. 
The calculation actually works for a general mapping stack $\mathbf{Map}(\Sigma,F)$ with a $d$-oriented stack $(\Sigma,[\Sigma])$ and an $n$-shifted symplectic stack $(F,\omega)$ as in Theorem \ref{thm-ptvv}. For a $f:\Sigma\to F$ of the mapping stack, the tangent complex at $f$ is 
$$
\mathbb{T}_f\mathbf{Map}(\Sigma,F)=\Gamma(\Sigma,f^*\mathbb{T}_F)\,.
$$
At $f$, the transgressed $(n-d)$-shifted symplectic form looks as follows: 
$$
\Gamma(\Sigma,f^*\mathbb{T}_F)^{\otimes 2}\overset{f^*\omega}{\longrightarrow}C^\infty(\Sigma)[n]\overset{[\Sigma]}{\longrightarrow}\mathbb{R}[n-d]\,.
$$
\begin{Exa}
Let $\Sigma$ be either the de Rham or Betti stack of a $d$-dimensional compact oriented manifold $X$, and let $F=BG$ for a Lie group $G$. 
Then $\Gamma(\Sigma,f^*\mathbb{T}_F)\cong H^*\big(X,ad(P_f)\big)[1]$, where $P_f$ is the flat $G$-bundle (or $G$-local system) on $\Sigma$ corresponding to the classifying map $f:X_{DR}\to BG$ (or $X_B\to BG$), $ad(P_f):=P_f\times_G\mathfrak g$ is the adjoint flat bundle
(or local system), and $H^*(X,-)$ means de Rham (or Betti) cohomology with coefficients. 

Setting $d=2$, we thus get a genuine linear symplectic pairing on the degree $0$ cohomology of $\Gamma(\Sigma,f^*\mathbb{T}_F)$, which is nothing but $H^1\big(X,ad(P_f)\big)$. This pairing is precisey the one we have seen in the Introduction (in the ``deformation theoretic approach'' part). 
\end{Exa}

\subsection{Transgression with boundary}

We have seen two systematic ways of constructing new shifted symplectic stacks out of old ones: 
\begin{itemize}
\item by doing derived intersection of Lagrangian morphisms. 
\item by transgression
\end{itemize}
We would like these two constructions to be compatible with each other. 
More precisely, we would like have the following property: if $M\cong M_+\coprod_N M_-$ is an oriented $d$-dimensional compact manifold obtained as the gluing of two manifolds sharing a common boundary $N\cong\partial M_+\cong\partial M_-$, and if $F$ is $n$-shifted symplectic, then both $\mathbf{Map}(M_B,F)$ and $\mathbf{Map}(M_{DR},F)$ can be obtained as derived Lagrangian intersections. 

More generally we will see that the transgression procedure produces a functor from the $\infty$-category of cobordisms to the $\infty$-category of shifted symplectic stacks and Lagrangian correspondences. 

\subsubsection{Boundary structures}

Let us summarize the structure we need in order to model the situation of a boundary inclusion in our framework. 
Let $\varphi:\Sigma\to \Upsilon$ be a morphism between derived stacks satisfying condition $1$ in paragraph \ref{integration}, and assume that $\Sigma$ is equipped with a $d$-shifted pseudo-orientation. 
\begin{Def}\label{def-bndstr}
A \textit{boundary structure} for $(\varphi,[\Sigma])$ is a homotopy $[\Upsilon]$ between $\varphi_*[\Sigma]:=[\Sigma]\circ\varphi^*$ and $0$ (as morphisms of cochain complexes $C^\infty(\Upsilon)\to \mathbb{R}[-d]$). 
\end{Def}
Below we list several examples. 
\begin{Exa}[Pseudo-orientations as boundary structures]\label{ex-pseudo}
Let $\emptyset$ be the initial stack, equipped with its canonical $d$-shifted orientation. 
A boundary structure on $\emptyset\to \Sigma$ is exactly the same as a $(d+1)$-shifted pseudo-orientation on $\Sigma$. 
This is very similar to the ``Symplectic is Lagrangian'' Example \ref{ex-symplag}. 
\end{Exa}
\begin{Exa}[Boundary inclusions as boundary structures]\label{ex-bndincl}
Let $N$ be a compact $(d+1)$-dimensional manifold with oriented boundary $M=\partial N$. Recall that both $M_B$ and $M_{DR}$ carry a $d$-shifted orientation. We further have that both maps $M_{DR}\to N_{DR}$ and $M_B\to N_B$, induced by the boundary inclusion $\iota:M\hookrightarrow N$, carry boundary structures. In the de Rham stack case this is simply Stokes formula $\int_{\partial M}\iota^*\omega=\int_Nd_{dR}\omega$: the homotopy is $\int_N$. In the Betti case, $[M_B]$ is given by the cap product with a fundamental cycle for $M$, and the homotopy $[N_B]$ is given by the cap-product with a compatible fundamental chain for $(N,M)$. 
\end{Exa}
\begin{Exa}\label{ex-obvious}
Let $M$ be a closed oriented surface, and recall from Example \ref{ex-curvature} that $B[\wedge^2T_M]$ is $1$-pseudo-oriented. 
The projection $\pi:\Sigma\to M$ carries an obvious boundary structure, as $\int_Mf=0$ if $f\in C^\infty(M)$. 
Indeed
\begin{itemize}
\item $\pi^*$ is just the inclusion $C^\infty(M)\hookrightarrow C^\infty(M)\oplus\Omega^2(M)[-1]$. 
\item the $1$-class $[\Sigma]$ is given by $\int_M$ (i.e.~it vanishes on functions and sends a $2$-form on $M$ to its integral). 
\item hence $\pi_*[\Sigma]=[\Sigma]\circ \pi^*$ necessarily vanishes (strictly). 
\end{itemize} 
\end{Exa}
\begin{Exa}\label{ex-Mconn}
Let $M$ be a closed oriented surface and consider again $\Sigma=B(\wedge^2T_M)$ with its $1$-shifted pseudo-oriented from Example \ref{ex-curvature}. We introduce another stack $M_{\nabla}$, being to $M_{DR}$ what connections are to flat connections: $M_{\nabla}$ is the quotient by the groupoid $\mathcal{G}\rightrightarrows M$ generated by the equivalence relation of being close at order $1$.\footnote{There is yet another description it terms of Lie algebroid on $M$, using the fact that every Lie algebroid leads to a formal thickening of $M$ (see \cite{GG,CG,Nuit}): 
\begin{itemize}
\item $M$ itself is associated with the trivial Lie algebroid $0$. 
\item $M_{DR}$ is associated with the Lie algebroid $T_M$. 
\item $M_{\nabla}$ is associated with the free Lie algebroid (see \cite{Kap}) generated by the anchored module $id:T_M\to T_M$. 
\item $B(\wedge^2T_M)$ is associated with the free Lie algebroid generated by the anchored module $0:\wedge^2 T_M\to T_M$. 
\end{itemize}
}
We have a morphism $\varphi:\Sigma\to \Upsilon:=M_{\nabla}$ that is roughly given by the map sending a connection to its curvature, which can be described in several ways: 
\begin{itemize}
\item There is a map that sends every infinitesimal $2$-simplex $\Delta^2_{inf}\to M$ to the element in $\mathcal{G}$ given by the sequence of equivalences $x\sim y\sim z\sim x$, where $x,y,z$ are the three vertices of the infinitesimal simplex, which induces a morphism on the quotients $B(\wedge^2TM)\to M_\nabla$. 
\item In term of Lie algebroids, we have a Lie algebroid morphism $Free(\wedge^2 T_M)\to Free (T_M)$ sending $u\wedge v$ to $uv-vu-[u,v]$, leading to a morphism $B(\wedge^2TM)\to M_{\nabla}$. 
\end{itemize}
At the level of functions, $C^\infty(\Upsilon)$ is the two-term complex $C^\infty(M)\overset{d_{dR}}{\longrightarrow}\Omega^1(M)$, $C^\infty(\Sigma)$ is the two term complex $C^\infty(M)\overset{0}{\longrightarrow}\Omega^2(M)$, and the morphism looks as follows: 
$$
\xymatrix{
C^\infty(\Upsilon) \ar[r]^{\varphi^*}& C^\infty(\Sigma) & \mathrm{degree} \\
C^\infty(M) \ar[r]^{id} \ar[d]_{d_{dR}} & C^\infty(M) \ar[d]_{0} & 0\\
\Omega^1(M) \ar[r]^{d_{dR}} & \Omega^2(M) & 1
}
$$
Hence $\varphi$ carries a boundary structure as the integral of an exact $2$-form on a closed manifold $M$ is zero. 
\end{Exa}

\noindent The following has been shown in \cite[Claim 2.7]{Cal15}: 
\begin{Prop}
Given a boundary structure as in Definiton \ref{def-bndstr} above, $\int_{[\Upsilon]}ev^*(-)$ provides a homotopy between the pull-back along $\varphi^*:\mathbf{Map}(\Upsilon,F)\to \mathbf{Map}(\Sigma,F)$ of $\int_{[\Sigma]}ev^*(-)$ and $0$. In particular, if $F$ is equipped with an $n$-shifted pre-symplectic structure then $\varphi^*$ carries an isotropic structure. 
\end{Prop}
We now provide several incarnations of this very general fact. 
\begin{Exa}[An infinite dimensional moment map for the moduli of connections]
Applying the above to the boundary structure on $M[\wedge^2T_M]\to M_{\nabla}$ from \ref{ex-Mconn} and to the $2$-shifted pre-symplectic structure on $F=BG$ ($G$ compact and simply connected) arising from $c\in S^2(\mathfrak g^*)^G$, we get an isotropic structure on 
the curvature morphism 
$$
\mathbf{Conn}_G(M):=\mathbf{Map}(M_{\nabla},BG)\to\mathbf{Map}\big(B(\wedge^2T_M),BG\big)=\big[\Omega^2(M,\mathfrak g)/C^\infty(M,\mathfrak g)\big]\,.
$$
\end{Exa}

If the boundary structure turns out to be non-degenerate in an appropriate sense, defining a \textit{relative $d$-orientation} (we refer to \cite[Definition 2.8]{Cal15} for the details) then, in complete analogy with Theorem \ref{thm-ptvv}, we have that the isotropic structure on $\varphi^*:\mathbf{Map}(\Upsilon,F)\to \mathbf{Map}(\Sigma,F)$ is a Lagrangian morphism\footnote{Whenever both mapping stacks are Artin, as usual. }. 
\begin{Exa}[Orientations as relative orientations]
Going back to Example \ref{ex-pseudo}, we have that a shifted relative $d$-orientation on $\emptyset\to \Sigma$ is exactly the same as a $(d+1)$-shifted orientation on $\Sigma$. Therefore if $(F,\omega)$ is an $n$-shifted symplectic stack then we have a Lagrangian stucture on $\mathbf{Map}(\Sigma,F)\to \mathbf{Map}(\emptyset,F)=*$, recovering Theorem \ref{thm-ptvv} from its relative analog. 
\end{Exa}
\begin{Exa}[Boundary inclusions as relative orientations]
The boundary structures from Example \ref{ex-bndincl} are always non-degenerate. Therefore, we obtain that whenever $G$ is a compact Lie group, then 
\begin{itemize}
\item $\mathbf{Loc}_G(M):=\mathbf{Map}(M_B,BG)$ is $(2-d)$-shifted symplectic and the restriction morphism $\mathbf{Loc}_G(N)\to\mathbf{Loc}_G(M)$ is Lagrangian. 
\item $\mathbf{Flat}_G(M):=\mathbf{Map}(M_{DR},BG)$ is $(2-d)$-shifted symplectic and the restriction morphism $\mathbf{Flat}_G(N)\to\mathbf{Flat}_G(M)$ is Lagrangian. 
\end{itemize}
\end{Exa}

\subsubsection{A gluing formula: transgression as a topological field theory}

Just like derived intersections of Lagrangian morphisms (resp.~isotropic morphisms) are shifted symplectic (resp.~shifted pre-symplectic), push-outs of relatively oriented morphisms (resp.~morphisms equipped with boundary structures) are oriented (resp.~pseudo-oriented). We refer to \cite[\S4.2.1]{Cal15} for the details. 

There is actually an $\infty$-category of oriented stacks with morphisms being relatively oriented cospans 
\[
  \xymatrix{
     \Sigma_1 \ar[dr] && \Sigma_2 \ar[dl] \\
    & \Upsilon & 
  }
\]
which works in a way similar (but dual) to our derived/$\infty$-categorical variant of Weinstein's symplectic category from \S\ref{weinstein-cat}, composition being given by push-out: 
\[
  \xymatrix{
    \Sigma_1 \ar[dr] && \Sigma_2 \ar[dl] \ar[dr] && \Sigma_3 \ar[dl] \\
    & \Upsilon_{12} \ar[dr] && \Upsilon_{23} \ar[dl] & \\
    && \Upsilon_{12}\underset{\Sigma_2}{\coprod}\Upsilon_{23} &&
  }
\]

We provide several examples in order to give the reader an intuition of what is going on. 
We restrict our attention to situations where $\Sigma_1=\Sigma_3=\emptyset$: then the composition is a relatively oriented cospan from $\emptyset$ to itself, and thus it is a $(d+1)$-oriented stack. 

\begin{Exa}[Betti and de Rham gluings]\label{ex-BDRgluing}
Let $N\cong N_+\coprod_M N_-$ be an oriented $(d+1)$-dimensional compact manifold obtained as the gluing of two oriented manifolds sharing a common boundary $M\cong\partial N_+\cong\partial N_-$. 

We therefore have $d$-oriented stacks $\Sigma=M_B$ and relative $d$-orientations on $\Sigma\to \Upsilon_{\pm}:=(N_{\pm})_B$. 
One can show that we have an equivalence of $(d+1)$-oriented stacks 
$$
\Upsilon_+\coprod_\Sigma\Upsilon_-\cong N_B\,.
$$
The same result holds as well for de Rham stacks. 
\end{Exa}

\begin{Exa}[Flat connections as connections with zero curvature]\label{ex-flatandco}
Let $M$ be an oriented compact surface. We have a push-out square of derived stacks\footnote{It comes from the following push-out square of Lie algebroids/groupoids over $M$: 
$$
\xymatrix{
Free(0:\wedge^2T_M\to T_M) \ar[r]\ar[d] & Free(id:T_M\to T_M) \ar[d] \\
0 \ar[r] & T_M
}
\qquad
\xymatrix{
(\wedge^2T_M,+) \ar[r]\ar[d] & \mathcal{G} \ar[d] \\
M \ar[r] & \widehat{M\times M}
}
$$
} 
$$
\xymatrix{
B(\wedge^2T_M)\ar[r]\ar[d] & M_{\nabla} \ar[d] \\
M \ar[r] & M_{DR}
}
$$
Recall that: 
\begin{itemize}
\item $B(\wedge^2T_M)$ is $1$-pseudo-oriented (Example \ref{ex-curvature}). 
\item $B(\wedge^2T_M)\to M_{\nabla}$ carries a boundary structure (Example \ref{ex-Mconn}). 
\item $B(\wedge^2T_M)\to M$ carries a boundary structure (Example \ref{ex-obvious}). 
\end{itemize}
One therefore gets that the push-out $M_{DR}$ is $2$-pseudo-oriented. 
We claim that the $2$-shifted pseudo-orientation coincides with the $2$-orientation on $M_{DR}$ from Example \ref{ex-DR}. 
\end{Exa}

One can show that the trangression procedure sends compositions of oriented cospans (resp.~of cospans with boundary structures) to compositions of Lagrangian (resp.~isotropic) correspondences. Restricting it to de Rham or Betti stacks we then get a $3d$-oriented topologicial field theory with values in our derived/$\infty$-variant of Weinstein's symplectic category. This $3d$ TFT can even be shown to be fully extended (see \cite{CHS}). 

\medskip

We again provide several examples. 

\begin{Exa}
We will apply the transgression procedure to the situation of Example \ref{ex-BDRgluing} with target $2$-shifted symplectic stack $BG$, with $G$ being compact. We thus get that the $(n-d-1)$-shifted symplectic stacks $\mathbf{Loc}_G(N)$ and $\mathbf{Flat}_G(N)$ that can be obtained as derived Lagrangian intersections: 
\[
\xymatrix{
\mathbf{Loc}_G(N) \ar[r]\ar[d] & \mathbf{Loc}_G(N_-) \ar[d]^{Lag} \\
\mathbf{Loc}_G(N_+) \ar[r]_{Lag} & \mathbf{Loc}_G(M)
}
\quad~~~~\quad
\xymatrix{
\mathbf{Flat}_G(N) \ar[r]\ar[d] & \mathbf{Flat}_G(N_-) \ar[d]^{Lag} \\
\mathbf{Flat}_G(N_+) \ar[r]_{Lag} & \mathbf{Flat}_G(M)
}
\]
\end{Exa}
\begin{Exa}
Let us specialize the above Example to the case when $N$ is a surface, $M=S^1$ and $N_-=D^2$ is a disk. We get that $\mathbf{Loc}_G(M)\cong[G/G]$ with its $1$-shifted symplectic structure and $\mathbf{Loc}_G(N_-)\cong[*/G]$. Hence we get back that the derived stack 
$\mathbf{Loc}_G(N)$ of $G$-local systems on the surface $N$ can be obtained as a quasi-Hamiltonian reduction of the derived stack of $G$-local systems on the same surface with a disk removed. I.e.~we have an equivalence of $0$-shifted symplectic stacks
$\mathbf{Loc}_G(N)\cong \mathbf{Loc}_G(N\backslash D^2)\times_{[G/G]}[*/G]$:  
\[
\xymatrix{
\mathbf{Loc}_G(N) \ar[r]\ar[d] & [*/G] \ar[d]^{Lag} \\
\mathbf{Loc}_G(N\backslash D^2) \ar[r]_-{Lag} & [G/G]
}
\]
If we also further assume that $G$ is connected then we have that $\mathbf{Flat}_G(M)\cong\big[\Omega^1(S^1,\mathfrak g)/C^\infty(S^1,G)\big]$ and $\mathbf{Flat}_G(N_-)\cong[*/G]$. Hence we get back the fact that the derived stack 
$\mathbf{Flat}_G(N)$ of flat $G$-bundles on the surface $N$ can be obtained as a kind of infinite dimensional Hamiltonian reduction of the derived stack of flat $G$-bundles on the same surface with a disk removed: we have a derived Lagrangian intersection 
\[
\xymatrix{
\mathbf{Flat}_G(N) \ar[r]\ar[d] & [*/G] \ar[d]^{Lag} \\
\mathbf{Flat}_G(N\backslash D^2) \ar[r]_-{Lag} & \big[\Omega^1(S^1,\mathfrak g)/C^\infty(S^1,G)\big]
}
\]
\end{Exa}

\begin{Exa}
We now apply the transgression procedure to the situation of Example \ref{ex-flatandco} with target $2$-shifted symplectic stack $BG$ again, $G$ being compact and simply connected for the sake of simplicity. We get a $1$-shifted pre-symplectic stack 
$\big[\Omega^2(M,\mathfrak g)/C^\infty(M,G)\big]$ together with isotropic morphisms to it from $\mathbf{Conn}_G(M)$ and $\mathbf{Bun}_G(M)=[*/C^\infty(M,G)]$. Their derived isotropic intersection is then equivalent, as a $0$-shifted pre-symplectic stack, to the $0$-shifted symplectic stack $\mathbf{Flat}_G(M)$: 
\[
\xymatrix{
\mathbf{Flat}_G(M) \ar[r]\ar[d] & \big[*/C^\infty(M,G)\big] \ar[d]^{isot.} \\
\mathbf{Conn}_G(M) \ar[r]_-{isot.} & \big[\Omega^2(M,\mathfrak g)/C^\infty(M,G)\big]
}
\]

\end{Exa}

\section*{Conclusion}
\addcontentsline{toc}{section}{Conclusion}

In this survey of derived symplectic geometry we have shown how one can use derived geometry in order to unify several descriptions of the symplectic structure on the reduced phase space of classical Chern--Simons theory, that is to say the moduli space of flat $G$-bundles, let's say for a simply connected compact Lie group, on a closed oriented surface $M$. 
Let us summarize what we have seen so far: 
\begin{itemize}
\item The derived moduli space $\mathbf{Flat}_G(M)$ is \textit{defined} as the derived mapping stack $\mathbf{Map}(M_{DR},BG)$. It naturally gets a $0$-shifted symplectic structure by transgression, using the $2$-shifted symplectic structure on $BG$ and the $2$-shifted orientation on $M_{DR}$. 
\item The canonical map $M_{DR}\to M_B$ of $2$-shifted oriented stacks induces an equivalence 
$\mathbf{Loc}_G(M)\tilde\longrightarrow\mathbf{Flat}_G(M)$ of $0$-shifted symplectic stacks, where 
$\mathbf{Loc}_G(M):=\mathbf{Map}(M_B,BG)$ is the derived stack of $G$-local systems on $M$. 
\item Locally at a point, the $0$-shifted symplectic structure is given by a combination of the non-degenerate pairing on $\mathfrak g$ with Poincar\'e duality. This can be seen as a derived extension of the standard fact that $H^1(M,\mathfrak g)$ is a symplectic vector space. 
\item The $0$-shifted symplectic structure on $\mathbf{Loc}_G(M)$ can be computed through the quasi-Hamiltonian reduction 
procedure of \cite{AMM}, being understood as a derived Lagrangian intersection. 
This is derived from a very general compatibility of the transgression procedure with all kinds of ``gluings''. 
\item Using a rather surprising gluing data (exhibiting $M_{DR}$ as a kind of push-out of $M_{\nabla}$) we can also compute the $0$-shifted symplectic structure on $\mathbf{Flat}_G(M)$ as an infinite dimensional reduction procedure applied to the derived moduli stack of all $G$-connections. This reduction procedure is here understood as a derived isotropic intersection. 
\item Whenever $M$ bounds an oriented $3d$-manifold $N$, we get a Lagrangian structure on the ``boundary condition map'' $\mathbf{Flat}_G(N)\to\mathbf{Flat}_G(M)$. 
\item In the case when $N$ is without boundary (i.e.~$M=\emptyset$) we can still make sense of what it means for the derived stack $\mathbf{Flat}_G(N)$ of boundary conditions to be Lagrangian in $*=\mathbf{Flat}_G(\emptyset)$: it means that $\mathbf{Flat}_G(N)$ is $(-1)$-shifted symplectic. 
\end{itemize}

\subsubsection*{What do we gain with shifted symplectic structures?}

As we have seen, derived geometry provide new tools, at the disposal of symplectic geometers, as well as a rather large framework in which several constructions find a nice interpretation, and do make sense in a wider context. 

We would like to point out several differences between ordinary symplectic geometry and its derived counterpart, that we may put in three different boxes:  
\begin{itemize}
\item first, there are differences coming from the use of Higher Algebra: several properties (e.g.~being closed, or isotropic) became structures. This has very interesting and unexpected consequences, such as the derived intersection of two lagrangian carrying a shifted symplectic structure. 
\item then, there are differences that actually come from the very purpose of derived geometry: several constructions do make sense in derived geometry even for pathological situations (including bad quotients, and non-transverse intersections). This allows for instance to deal with symplectic reduction without any hypothesis on the action, or to define Weinstein's symplectic category. 
\item finally, there are less expected differences: we discover the existence of a new kind symplectic structures, that were somehow hidden in ordinary symplectic geometry: the $2$-shifted symplectic structure on $BG$, the $1$-shifted symplectic structure on $[\mathfrak{g}^*/G]$, or moment maps being lagrangian morphisms. 
\end{itemize}
There is one particular instance of a construction that reveals all three aspects that we have just listed: it is the transgression Theorem \ref{thm-ptvv} and its generalizations. Indeed, the proof makes a crucial use of structures rather than properties, in many cases the mapping spaces involved become tractable derived stacks (avoiding the use infinite dimensional geometry and the need to restrict on some smooth locus), and we get access to negatively shifted symplectic structures just from the $2$-shifted symplectic structure on $BG$. 

\medskip

Let us finally mention a very interesting consequence of being $(-1)$-shifted symplectic for $\mathbf{Flat}_G(N)$, $N$ being a closed oriented $3$-manifold. In \cite{BBBJ} it has been proved that $(-1)$-shifted derived stacks locally look like derived critical loci. In particular this implies that the corresponding \textit{non-derived} Artin stack $Flat_G(N)$ is a $d$-critical stack in the sense of \cite{Jo}. Also note that there are sufficient conditions for the existence of a global action functional for $(-1)$-shifted derived stacks \cite{Pantev}. For instance, it is expected that if a $(-1)$-shifted derived stack carries a Lagrangian foliation $\mathcal L$ then there exists a function $f$ on the formal quotient stack $[X/\mathcal L]$ and an \'etale symplectic morphism $X\to \mathbf{Crit}(f)$. This means that somehow, the action functional is hidden in the $(-1)$-shifted symplectic structures, which may be a more fundamental structure. 

\subsubsection*{Shifted Poisson structures}

In ordinary symplectic geometry, any symplectic manifold carries a corresponding Poisson structure. A natural question is then: \textit{are there shifted Poisson structures as well?} And the answer is \textit{yes}. Shifted Poisson structures have been introduced and studied in \cite{CPTVV,Pri}, where it is shown that non-degenerate shifted Poisson structures do coincide with shifted symplectic structures. Several observations are in order: 
\begin{itemize}
\item the theory of shifted Poisson structures is rather technical, involving complications due to the lack of functoriality of the tangent complex (as opposed to the cotangent complex). Simple facts like the bijection between non-degenerate Poisson bivectors and symplectic forms requires quite a lot of work in the derived setting.  
\item many constructions and structures from standard Poisson geometry (such as Poisson-Lie structures, $r$-matrices, dynamical $r$-matrices, the Feigin--Odesski algebra, $\dots$) fit very well in the realm of shifted Poisson structures, as it has been shown in Safronov's recent work \cite{Saf17}. 
\item $(-1)$-shifted Poisson structures are of particular interest for the classical BV--BRST formalism. In particular, using \cite{CPTVV,Pri} one can show that the function ring of $\mathbf{Flat}_G(N)$ carries a Poisson bracket of cohomological degree $1$: this bracket indeed lies at the heart of the BV--BRST formalism, and is the starting point of perturbative quantization for several classical field theories. More details about the relation between shifted Poisson geometry and the classical and quantum BRST formalism can be found in \cite{Safred}. 
\item there is an alternative definition of shifted Poisson structures as formal Lagrangian thickenings. The idea being that if $X$ is Poisson then the formal quotient $[X/\mathcal L]$ by its symplectic foliation $\mathcal L$ is $1$-shifted symplectic and the map $X\to[X/\mathcal L]$ is Lagrangian. At present it is still a conjecture that shifted Poisson structures are equivalent to formal Lagrangian thickenings. 
\end{itemize}

\subsubsection*{Quantization}

Having new kinds of (pre)symplectic and Poisson structures at our disposal, we may want to know if there is a notion of quantization for these, especially for examples coming from classical field theory. The theory is still at its infancy, but there has already been a few progresses: 
\begin{itemize}
\item deformation quantization of $n$-shifted Poisson structures has been discussed in \cite{CPTVV}, and makes sense for $n\geq-1$. It has been shown to exist whenever $n>0$ in \cite{CPTVV}. The cases $n=-1$ (BV quantization, or rather BD quantization after \cite{CoGw}) and $n=0$ (usual case) are the most difficult to deal with (we refer to Pridham's work \cite{Pri1,Pri2} for partial results). 
\item geometric quantization has been initiated in a paper of James Wallbridge \cite{Wa} (see also the end of \cite{Saf16}). 
\end{itemize}
In \cite{CPTVV,Pri} there is a notion of compatible pair of an $n$-shifted presymplectic structure and an $n$-shifted Poisson structure. It would be interesting to investigate a notion of compatible pair of quantizations (deformation and geometric). Indeed, even in the non-derived situation, comparing the two kinds of quantizations is a relevant question.  

\newcommand{\bysame}{\leavevmode\hbox to3em{\hrulefill}\,}

\end{document}